\documentclass[10pt, final]{ieeetran}

\usepackage{graphicx}
\usepackage[dvips]{epsfig}

\usepackage{amssymb,amsmath,amsfonts,enumerate}
\usepackage{cite}
\usepackage{caption,subcaption}
\usepackage{color}

\usepackage{hyperref}%
\hypersetup{colorlinks=true, linkcolor=blue, breaklinks=true, urlcolor=blue, citecolor=blue}

\def\la{\lambda}
\def\La{\Lambda}
\def\vp{\varphi}
\def\ve{\varepsilon}
\def\vk{\varkappa}

\def\g{\mathcal G}

\def\r{\mathbb R}

\def\be{\beta}

\def\ones{\mathbf 1}
\def\conv{\mathop{conv}\nolimits}
\def\dist{\mathop{dist}\nolimits}

\newcommand{\dfb}{\stackrel{\Delta}{=}}

\def\E{\mathcal E}
\def\V{\mathcal V}

\def\be{\begin{equation}}
\def\ee{\end{equation}}
\def\ben{\begin{equation*}}
\def\een{\end{equation*}}

\newtheorem{thm}{Theorem}
\newtheorem{lem}{Lemma}
\newtheorem{cor}{Corollary}
\newtheorem{defn}{Definition}
\newtheorem{rem}{Remark}
\newtheorem{assum}{Assumption}
\newtheorem{prop}{Proposition}

\title{Delay Robustness of Consensus Algorithms:\\ Beyond The Uniform Connectivity\\[2mm] \small\bf(Extended Version)}

\author{Anton V. Proskurnikov and Giuseppe Carlo Calafiore%
\thanks{The authors are with the Department of Electronics and Telecommunications, Politecnico di Torino, Turin, Italy.
Anton V. Proskurnikov is also with Institute for Problems of Mechanical Engineering, Russian Academy of Sciences,
St. Petersburg, Russia. %The work is supported by RFBR grant 20-01-00619.
Email: {\tt\small anton.p.1982@ieee.org, giuseppe.calafiore@polito.it}
}%
}

\begin{document}

\maketitle
\thispagestyle{empty}
\pagestyle{empty}

\begin{abstract}
Consensus of autonomous agents is a benchmark problem in multi-agent control. In this paper, we consider continuous-time averaging consensus policies (or Laplacian flows) and their
discrete-time counterparts over time-varying graphs in presence of unknown but bounded communication delays. It is known that consensus is established (no matter how large the delays are) if the graph
is periodically, or uniformly quasi-strongly connected (UQSC).
The UQSC condition is often believed to be the weakest sufficient condition under which consensus can be proved.
We show that the UQSC condition can actually be substantially relaxed and replaced by a condition that we call \emph{aperiodic} quasi-strong connectivity (AQSC),
which, in some sense, proves to be very close to the necessary condition of integral connectivity.
Furthermore, in some special situations such as undirected or type-symmetric graph,
we find a necessary and sufficient condition for consensus in presence of bounded delay; the relevant results have been previously proved only in the undelayed case.
The consensus criteria established in this paper generalize a number of results known in the literature.
%Furthermore, consensus can be proved for any feasible solution of the delay differential inequalities associated to the consensus algorithm. Such inequalities naturally arise in problems of containment control, %distributed optimization and models of social dynamics.
\end{abstract}

%%%%%%%%%%%5
\section{Introduction}

Consensus policies are prototypic distributed algorithms for multi-agent coordination~\cite{RenBeardBook,RenCaoBook}
inspired by regular ``intelligent'' behaviors of biological and physical systems~\cite{Vicsek,Strogatz:00,Couzin:2011}.
The most studied first-order consensus algorithms are based on the principle of iterative averaging; such algorithms were considered~\cite{French:1956,Abelson:1964,DeGroot} long before the recent
``boom'' in multi-agent control and networks. Averaging policies have found numerous applications
as exemplified by models of sociodynamics~\cite{ProTempo:2018} and algorithms for distributed computing~\cite{Tsitsiklis:86,ShiJohanssonHong:13,Nedic:10,Wang2019_ARC,FullmerMorse2018,Amelina:15}.

Consider a finite team of agents $\V$, each of which is associated with some value of interest $x_i\in\r,\,i\in\V$ representing for instance attitudes~\cite{Abelson:1964} or the value of some physical characteristics~\cite{Murray:07}.
In the discrete-time case, the agents simultaneously update their values to the \emph{average} of their  own value and and the other agents' values:
\be\label{eq.conse0}
\begin{aligned}
&x_i(t+1)=\sum_{j\in\V} a_{ij}(t)x_j(t)=\\&=x_i(t)+\sum_{j\ne i}a_{ij}(t)[x_j(t)-x_i(t)],
\forall i\in\V,t=0,1,\ldots
\end{aligned}
\ee
where $(a_{ij}(t))\in\r^{\V\times\V}$ is a \emph{stochastic} matrix (all entries are nonnegative, each row sums to $1$). The continuous-time counterpart of~\eqref{eq.conse0}, called the \emph{Laplacian flow~\cite{BulloBook-Online}}) , is
\be\label{eq.conse1}
\dot x_i(t)=\sum_{j\ne i} a_{ij}(t)(x_j(t)-x_i(t))\quad\forall i\in\V,\,t\geq 0
\ee
where $(a_{ij}(t))\in\r^{\V\times\V}$ is a nonnegative matrix (not necessarily stochastic). 
Entry $a_{ij}(t)$ is interpreted as the weight of influence of agent $j$ on agent $i$ at time $t$: the larger 
weight is, the stronger is attraction of agent $i$'s value to agent $j$'s value.

The central question regarding dynamics~\eqref{eq.conse1} is establishing eventual \emph{consensus}, that is, convergence of all values $x_i(t)$ to the same value $\bar x=\lim_{t\to\infty}x_i(t)\,\forall i$ (which may depend on the initial condition).  More general behavior is ``partial'' consensus (group consensus, clustering)~\cite{Yu10,XiaCao:11}, that is, splitting of the agents into several groups that agree on different values. In the case of constant weights $A(t)\equiv A$, the consensus criterion is well-known~\cite{ChebotarevAgaev:2002,RenBeardBook} and reduces to the quasi-strong connectivity of the graph associated to matrix $A$ (equivalently, the graph has a directed spanning tree). Informally, some agent influences all other agents in the network, directly or indirectly. Without this ``weak'' connectivity, the values $x_i(t)$ converge, and their steady values are determined by the graph's spanning forest~\cite{AgaevChe:2014}. In the discrete-time case, consensus is equivalent to ergodicity of the backward infinite products $A(t)\ldots A(0)$~\cite{Seneta,LEIZAROWITZ1992189}; this line of research, related to inhomogeneous Markov chain theory, was developed in the recent works~\cite{TOURI20121477,Bolouki:2016,Touri:14}.

Finding criteria ensuring consensus in the case of a general time-varying matrix $(a_{ij}(t))$ is a difficult problem whose complete solution is still elusive. A well-known \emph{necessary} condition for consensus is the so-called \emph{integral} (essential, persistent) connectivity~\cite{TsiTsi:13,MatvPro:2013,ShiJohansson:13-1}. Namely, the ``persistent'' arcs corresponding to such pairs of agents that
\be\label{eq.int-aij}
\sum_{t=0}^{\infty}a_{ij}(t)=\infty\quad\text{or}\quad \int_0^{\infty}a_{ij}(t)dt=\infty
\ee
(depending on whether the time is discrete or continuous),
should constitute a quasi-strongly connected graph. This condition, however, is far from being sufficient. Simple counter-examples~\cite{Moro:05,MartinHendrickx:2016} show that the agents may fail to reach consensus (and, moreover, vector $x$ may fail to converge) even if the graph of persistent interactions is complete. This effect is caused by ``imbalance'' between influence weights: some pairs of agents interact more intensively than others. In the discrete-time case, alternative necessary conditions inspired by theory of inhomogeneous Markov chains were proposed in~\cite{TOURI20121477,Bolouki:2016};
the simplest of these conditions termed~\emph{infinite flow}, in fact, follows from persistent connectivity, whereas the ``absolute''~\cite{TOURI20121477} and ``jet''~\cite{Bolouki:2016} infinite flow conditions
require that the connectivity persists under special transformations of the matrix $A(t)$, being thus highly non-trivial conditions of time-varying \emph{aperiodicity}. The verification of these properties is, however, a self-standing non-trivial problem, which has been solved only in special cases.

Sufficient criteria for consensus can be divided into two groups. Conditions of the first type require the periodic, or \emph{uniform} quasi-strong connectivity~\cite{Moro:04,RenBeardBook,LinFrancis:07}:
two numbers $T,\ve>0$ should exist such that the unions of the interacting graphs over each interval $[t,t+T],\,t\geq 0$ are quasi-strongly connected and this connectivity property persists if one removes ``lightweight'' arcs whose weights are less than $\ve$. Often, additional requirements like positive dwell-time between subsequent switchings~\cite{Antonis} are added. The uniform connectivity is not necessary for consensus in the usual sense and implies, in fact, much stronger properties of consensus with uniform~\cite{LinFrancis:07} or exponential convergence~\cite{BarabanovOrtega2018} and consensus robustness against unknown disturbances~\cite{ShiJohansson:13}.

Consensus criteria of the second kind ensure consensus in presence of the integral connectivity and some conditions preventing the imbalance of couplings. The simplest condition of this type is the coupling symmetry $a_{ij}(t)=a_{ji}(t)$~\cite{CaoZheng:11}, which condition can be in fact relaxed to weight-balance, type-symmetry or cut-balance conditions~\cite{TsiTsi:13}. All of these conditions guarantee \emph{reciprocity} of interactions: if some group of agents $S\subset\{1,\ldots,n\}$ influences the remaining agents from $S^c=\{1,\ldots,n\}\setminus S$, then agents from $S^c$ also influence agents from $S$, moreover, the mutual influences of groups $S$ and $S^c$ are commensurate. In the discrete-time case, generalizations of cut-balance conditions were proposed in~\cite{Touri:14,Bolouki2015,Bolouki:2016}.
Criteria for consensus under \emph{non-instantaneous} reciprocity conditions were obtained in~\cite{MartinHendrickx:2016,MartinGirard:2013,XiaShiCao2019,ProCalaCao:2020}.

An important question regarding consensus algorithms is robustness against \emph{communication delays}. Such delays naturally arise in the situation where the agents have direct access to their own values, whereas the neighbors' values are subject to non-negligible time lags. Delays of this type are inevitable in networks spread over large distances (e.g., where the agents communicate via Internet), but also arise in many physical models~\cite{StrogatzDelay1:03,Antonis,DahmsSchoell:2012}. Consensus criteria of the first type (based on the uniform connectivity) in fact ensure consensus robustness against arbitrary bounded delays~\cite{Moro:04,Bliman:06,Antonis,Muenz:11}.  However, delay robustness without uniform connectivity of the graph has remained an open problem; the existing results are mostly limited to the discrete-time case and impose restrictive conditions on the matrix $A(t)$, e.g., the uniform positivity of its non-zero entries~\cite{Blondel:05,ProCalaCao:2020}. There is a judgement that, dealing with delayed consensus algorithms, the ``uniform quasi-strong connectivity is in fact the weakest assumption on the graph connectivity such that consensus is guaranteed for arbitrary initial conditions,''~\cite{Muenz:11}.

In this paper, we extend the theory of delayed averaging algorithms in two directions. \textbf{First,} we show that in the case of a general time-varying directed graph, the uniform connectivity can be substantially relaxed to a condition termed below \emph{aperiodic} connectivity, which is very close to the necessary condition of persistent connectivity. This result, in fact, is of interest even for the undelayed case and, as it will be shown, generalizes some consensus conditions available in the literature~\cite{Muenz:11,ShiJohansson:13-1,XiaShiCao2019}.  \textbf{Second}, we extend the reciprocity-based consensus criteria established in~\cite{TsiTsi:13,MatvPro:2013,MartinHendrickx:2016} that ensure consensus when the graph is type-symmetric (in the ``non-instantaneous'' sense) and persistently connected.

In the undelayed case the discrete-time equation~\eqref{eq.conse0} is often considered to be simpler than the continuous-time dynamics~\eqref{eq.conse1}; many results for continuous-time algorithms are derived via discretization~\cite{RenBeardBook,Bolouki:2016,MartinHendrickx:2016}. In the case of delayed communication, the opposite is the case: the discrete-time averaging algorithm appears to be a special case of the infinite-dimensional continuous-time consensus dynamics. Using this fact, we develop a unified theory of continuous-time and discrete-time consensus algorithms.

The paper is organized as follows. Section~\ref{sec.prelim} introduces preliminary concepts and notation.
Section~\ref{sec.setup} provides the problem setup (delay-robust consensus in averaging algorithms);
we also discuss known necessary condition of consensus (persistent connectivity) there.
Section~\ref{sec.aqsc} offers the first sufficient condition for consensus, applicable to a general directed graph and
generalizing the commonly used UQSC condition.
Another sufficient condition, assuming the non-instantaneous type-symmetry, is introduced in Section~\ref{sec.symm}.
Section~\ref{sec.ext} discusses some further extensions of our consensus criteria, in particular, we show that
our criteria ensure in fact robustness of consensus not only against delays, but also against $L_1$-summable and some other
vanishing disturbances (Subsect.~\ref{subsec.robust}). We show that the theory developed in this paper can be extended,
without major differences, to nonlinear averaging algorithms with delays (Subsect.~\ref{subsec.nlin}). Finally,
we demonstrate in Subsect.~\ref{subsec.contain} that the theory of consensus algorithms, in fact, allows to examine more general algorithms
of multi-agent coordination referred to as the containment control~\cite{RenCaoBook} and target aggregation~\cite{ShiHong:09}
policies: the goal of such algorithms is to gather mobile agents in a polytope or some other predefined set.
Section~\ref{sec.proof} collects the technical proofs of the main results. Section~\ref{sec.concl} concludes the paper.

\section{Preliminaries}\label{sec.prelim}

Throughout the text, symbol $\dfb$ should be read as ``defined as''. For integers $m\leq n$, let $[m:n]\dfb\{m,m+1,\ldots,n\}$.

%\subsection{Vectors, matrices, graphs}

Given a finite set of indices $\V$, we use $\r^{\V}$ to denote the set of vectors $x=(x_i)_{i\in\V}$, where $x_i\in\mathbb{R}$.
For such a vector, $\min x\dfb\min_i x_i$ and $\max x\dfb\max_ix_i$. As usual, $\|x\|_{\infty}\dfb\max_i|x_i|$.
For two vectors $x,y\in\r^{\V}$, we write $x\leq y$ if $x_i\leq y_i\,\forall i$. We use $\r^{\V\times\V}$ to denote the set
of matrices $A=(a_{ij})_{i,j\in\V}$, where $a_{ij}\in\r\,\forall i,j\in\V$.

The vectors of standard coordinate basis in $\r^{\V}$ are denoted by $\mathbf{e}^i\dfb(\delta_j^i)_{j\in\V}$, where
$\delta_i^i\dfb 1\,\forall i$ and $\delta_j^i\dfb 0\,\forall j\ne i$.
Let $\ones_{\V}\dfb\sum_{i\in\V}\mathbf{e}^i$ denote the vector of ones and $I_{\V}\dfb (\delta^i_j)_{i,j\in\V}$ be the identity matrix; the subscript $\V$ will be omitted when this does not lead to confusion. A nonnegative matrix $A\in\r^{\V\times\V}$ is \emph{stochastic} if $A\ones=\ones$ and \emph{substochastic} if $A\ones\leq\ones$.

 A (directed) graph is a pair $\g=(\V,\E)$, where $\V$ is a finite set of \emph{nodes} and $\E\subseteq\V\times\V$ is the set of \emph{arcs}.  A \emph{walk} from node $i\in V$ to node $j\in V$ is a sequence of arcs $(v_0,v_1)$, $(v_1,v_2)$,\ldots,$(v_{n-1},v_n)$ starting at $i_0=i$ and ending at $i_n=j$. A graph is \emph{strongly connected} if every two nodes are connected by a walk and \emph{quasi-strongly connected} (has a directed spanning tree~\cite{RenBeardBook}, rooted~\cite{CaoMorse:08}) if some node (a \emph{root}) is connected to all other nodes by walks.

 A \emph{weighted} graph, determined by  a weight matrix $A\in\r^{\V\times\V}$,   is a triple $\g[A]\dfb(\V,\E,A)$, where $A\in\r^{\V\times\V}$ is a \emph{nonnegative} matrix
 that is compatible with graph $(\V,\E)$, that is\footnote{Following the convention adopted in multi-agent systems~\cite{RenBeardBook}, influence of agent $j$ on agent $i$ is represented by arc $(j,i)$ rather than $(i,j)$.}, $\E\dfb\{(j,i)\in\V\times\V:a_{ij}>0\}$. Given $\ve>0$, denote
 \[
 A^{[\ve]}\dfb(a_{ij}^{[\ve]}),\quad a_{ij}^{[\ve]}=
 \begin{cases}
 a_{ij},\quad a_{ij}\geq\ve,\\
 0,\quad a_{ij}<\ve.
 \end{cases}
 \]
 Graph $\g[A^{[\ve]}]$ is called the \emph{$\ve$-skeleton} of graph $\g[A]$ (this graph is obtained from $\g[A]$ by removing ``lightweight'' arcs of weight $<\ve$).
  We call a graph \emph{(quasi)-strongly $\ve$-connected} if its $\ve$-skeleton is (quasi-)strongly connected.

\section{Problem setup. Necessary conditions.}\label{sec.setup}

In this paper, we are interested in the delayed counterparts of the systems~\eqref{eq.conse0} and~\eqref{eq.conse1}: the continuous-time equation
\be\label{eq.conse1d}
\dot x_i(t)=\sum\nolimits_{j\ne i}a_{ij}(t)(\hat x_j^i(t)-x_i(t)),\,i\in\V,
\ee
and its discrete-time counterpart
\be\label{eq.conse0d}
x_i(t+1)=x_i(t)+\sum\nolimits_{j\ne i}a_{ij}(t)(\hat x_j^i(t)-x_i(t)),\,i\in\V,
\ee
Unlike the algorithms~\eqref{eq.conse0d} and~\eqref{eq.conse1d}, at time agent $i$ receives a \emph{retarded}\footnote{To simplify matters, we consider only discrete delays, the results, however, can be
extended to distributed bounded delays~\cite{Muenz2009} without essential changes. Consensus, in fact, can also be robust to unbounded delays~\cite{LiuLuChen:2010,XuLiuFeng:2019}, however, such extensions
are beyond the scope of this paper.} value of agent $j\ne i$ denoted by
\[
\hat x_j^i(t)\dfb x_j(t-{h}_{ij}(t)),\quad h_{ij}(t)\in [0,\bar h].
\]
The coefficients obey the standard assumption, which is henceforth \textbf{always} supposed to hold:
\begin{assum}\label{asm.coeff}
 The weights $a_{ij}(t)\geq 0$ and  delays $h_{ij}(t)\in [0,\bar h]$ satisfy the following conditions:
 %are defined for all $t\in[0,\infty)$. Also,
 \begin{itemize}
 \item in continuous-time algorithm~\eqref{eq.conse1d}, $a_{ij}$ are locally $L_1$-summable and $h_{ij}$ are measurable on $[0,\infty)$;
 \item in the discrete-time case, $\sum_{j\in\V}a_{ij}(t)=1\,\forall i\in\V$ and $h_{ij}(t)$ are integer for all $t=0,1,\ldots$
 \end{itemize}
\end{assum}

Notice that we formally confine ourselves to linear consensus algorithms with scalar values $x_i(t)\in\r$, however, the extensions to a broad class of nonlinear algorithms in the general vector
spaces prove to be straightforward (Section~\ref{sec.ext}).

In the definition of consensus, it is convenient to consider algorithms~\eqref{eq.conse1d} and~\eqref{eq.conse0d} with the starting time $t_*\geq 0$. The solution to~\eqref{eq.conse0d}  is uniquely defined
for $t>t_*$ by the initial condition  $x(t_*),x(t_*-1),\ldots,x(t_*-\bar h)$.  The solution to~\eqref{eq.conse1d}, which will always supposed to be right-continuous at $t_*$, is uniquely determined on $[t_*,\infty)$
by the initial condition
\be\label{eq.initial}
x(t_*)=x_*\in\r^n,\; x(t+s)=\vp(s)\,\forall s\in[-\bar h,0)
\ee
where $\vp\in L_{\infty}([-\bar h,0]\to\r^{\V})$ is some known function~\cite{Hale}.

\subsection*{Problem setup: delay-robust consensus}

We start with the definition of consensus.

\begin{defn} The algorithm~\eqref{eq.conse1d} (or~\eqref{eq.conse0d}) establishes consensus in the group of agents $\V'\subseteq\V$ if
for any initial condition~\eqref{eq.initial} (or its discrete-time counterpart) and
any agent $i\in\V'$ there exists the ultimate value
\be\label{eq.barx-i}
\bar x_i\dfb\lim_{t\to\infty} x_i(t),
\ee
and these ultimate values are coincident: $\bar x_i=\bar x_j\,\forall i,j\in\V'$. If this holds for $\V'=\V$, then \emph{global} consensus is established.
\end{defn}

Notice that the initial condition also includes the \emph{starting time} $t_*\geq 0$; in particular, our definition of consensus excludes non-generic situations where consensus
in the discrete-time algorithm is reached in finite time due to a specific structure of matrices $A(t)$~\cite{Hendrickx:2015}, after which communication between the agents can be completely stopped.
We are primarily interested in the criteria for global consensus, however, conditions for solution's convergence and consensus in subgroups will also be obtained as a byproduct (see Theorems~\ref{thm.symm-c}
and~\ref{thm.symm-d}).  Observe that due to time-varying weights and delays, the existence of limits~\eqref{eq.barx-i} in the time-varying system is a non-trivial self-standing problem.

Since delays are usually caused by imperfect communication channels, we suppose that they are uncertain, defining  thus the problem of \emph{delay-robust} consensus as follows.

\textbf{Problem.} Find conditions on the time-varying matrix $A(\cdot)$ ensuring that the algorithm~\eqref{eq.conse0d} or~\eqref{eq.conse1d} establishes global consensus \emph{for all} possible delays $h_{ij}(t)\in [0,\bar h]$ (in the continuous-time case, $h_{ij}$ are Lebesgue measurable).

As it will be discussed in Section~\ref{sec.ext}, our conditions in fact ensure also the robustness of consensus against unknown \emph{disturbances}, extending the results from~\cite{ShiJohansson:13} to the case
of delayed consensus algorithms.

Introducing the maximal and minimal values of the agents over the time window $[t-\bar h,t]$ as follows
\be\label{eq.la}
\la(t)\dfb\inf\limits_{t-\bar h\leq s\leq t}\min x(s),\quad\La(t)\dfb\sup\limits_{t-\bar h\leq s\leq t}\max x(s),
\ee
it can be shown (see Lemma~\ref{prop.bound}) that $\la$ and $\La$ are monotonically non-decreasing and non-increasing
respectively. Hence, global consensus can be alternatively formulated as
\be\label{eq.conse-diameter}
\La(t)-\la(t)\xrightarrow[t\to\infty]{} 0.
\ee

\subsection*{Necessary conditions for consensus}

Since in the undelayed case $\bar h=0$ the algorithms~\eqref{eq.conse1d} and~\eqref{eq.conse0d} reduce, respectively, to conventional averaging dynamics~\eqref{eq.conse1} and~\eqref{eq.conse0},
the necessary conditions for consensus without delays are \emph{a fortiori} necessary for delay-robust consensus. These conditions are formulated in terms of the  \emph{persistent} graph.
\begin{defn}
Agent $j$ persistently interacts with agent $i$ if one of the relations~\eqref{eq.int-aij} holds (corresponding to the discrete or continuous time). Denoting the set of such pairs $(j,i)$ by $\E_{\infty}$,
graph $\g_{\infty}=(\V,\E_{\infty})$ is said to be the graph of persistent interactions, or the  \emph{persistent graph} of the algorithm.
\end{defn}
\begin{lem}\label{lem.necess}~\cite{ShiJohansson:13-1}\footnote{Formally, consensus criteria in~\cite{ShiJohansson:13-1} are primarily focused on the conditions for global consensus, but
the necessary conditions for consensus between two agents are established as a part of their proofs.}
If algorithm~\eqref{eq.conse0} or~\eqref{eq.conse1} establishes consensus between two agents $i,j\in\V$ (that is, $\bar x_i=\bar x_j$), then at least one statements (a)-(c) is valid:
a) $(j,i)\in\E_{\infty}$; b) $(i,j)\in\E_{\infty}$; c) there exists some agent $k\ne i,j$ such that $(i,k),(j,k)\in\E_{\infty}$.
In particular, if the \emph{global} consensus is established, then graph $\g_{\infty}$ is quasi-strongly connected.
\end{lem}

The persistent graph's connectivity admits an equivalent reformulation, considering \emph{unions} of graphs $\g[A(t)]$. Let
\be\label{eq.unions}
\begin{aligned}
A_{t_1}^{t_2}\dfb\int_{t_1}^{t_2}A(s)\,ds\quad \text{for algorithm~\eqref{eq.conse1d}},\\
A_{t_1}^{t_2}\dfb\sum\nolimits_{t=t_1}^{t_2}A(s)\,ds\quad \text{for algorithm~\eqref{eq.conse0d}},
\end{aligned}
\ee
(in both situations, $0\leq t_1\leq t_2\leq\infty$, in the second case, $t_1,t_2$ are integer), we call the respective graph $\g[A_{t_1}^{t_2}]$ the \emph{union} of graphs over interval $[t_1,t_2]$ (respectively, $[t_1:t_2]$).
\begin{lem}\label{prop.necess}
The persistent graph $\g_{\infty}$ is quasi-strongly connected  if and only a constant $\ve>0$ and increasing sequence $t_p\xrightarrow[p\to\infty]{}\infty$ exists ($p=0,1,2,\ldots$) such that all unions of the graph $\g[A_{t_p}^{t_{p+1}}]$ are quasi-strongly $\ve$-connected. If both statements are valid, then $\ve$ can be chosen arbitrarily.
\end{lem}
\begin{IEEEproof}
Denoting $\g_p\dfb\g[A_{t_p}^{t_{p+1}}],\,p\geq 0$, let $\g_p^{[\ve]}$ stand for the $\ve$-skeleton of $\g_p$. To prove the \textbf{``if''} part, assume that every graph $\g_p^{[\ve]}$ are quasi-strongly connected and therefore contain directed (out-branching) spanning tree. Since the number of trees is finite, some of them belongs to an infinite subsequence of graphs $\g_{p_k}$, where $p_k\to\infty$, being thus contained by $\g_{\infty}$. To prove the \textbf{``only if''} part, consider the sequence $t_p$ constructed as follows: $t_0=0$ and for each $p\geq 0$ one has
\[
t_{p+1}=\inf\{t>t_p:\int_{t_p}^{t_{p+1}}a_{ij}(t)dt\geq\ve\quad\forall (j,i)\in\g_{\infty}\}
\]
(in the discrete-time case, the integral is replaced by the sum). By construction, each graph $\g_p^{[\ve]}$ contains $\g_{\infty}$, being thus quasi-strongly connected.
\end{IEEEproof}

In the next section, we formulate a sufficient condition for consensus which is close to the necessary condition from Lemma~\ref{prop.necess} yet imposes a restriction
on the matrices $A_{t_p}^{t_{p+1}}$.

\subsection*{An extra assumption: a stronger form of local $L_1$-summability}

In the continuous-time case, our consensus criteria will rely on an additional condition which cannot be omitted in the case of general time-varying delays (see Appendix~\ref{sec.app-a}):
\begin{assum}\label{asm.int-bound}
For any finite number $D>0$, one has
\be\label{eq.int-bound}
\mu_D\dfb\sup_{\substack{t\geq 0\\i,j\in\V}}\int\nolimits_t^{t+D}a_{ij}(s)\,ds<\infty.
\ee
\end{assum}

Also, it suffices to check~\eqref{eq.int-bound} for some arbitrary $D>0$.
Assumption~\ref{asm.int-bound}, obviously, holds if $a_{ij}$ are bounded functions and is often assumed to hold in consensus criteria~\cite{ShiJohansson:13}, even though it can often be discarded
in the undelayed case (see, e.g., the last statement of Theorem~\ref{thm.repeated-cont} below).

\section{The first sufficient condition:\\ Aperiodic quasi-strong connectivity}\label{sec.aqsc}

In this section, we establish the first sufficient condition for delay-robust consensus, based on the property of \emph{aperiodic} quasi-strong connectivity. We start with definitions.

\begin{defn}\label{def.tp-bound}
Consider an increasing sequence $\mathfrak{t}=(t_p)_{p=0}^{\infty}$, such that $t_p\to\infty$, $t_p\geq 0$ and, in the discrete-time case, $t_p$ are integer. The matrix function $A(\cdot)$ is \emph{$\mathfrak{t}$-bounded} if
\be\label{eq.ell}
\ell=\ell(A(\cdot),\mathfrak{t})=\sup_{\substack{p=0,1,\ldots\\i\ne j}}\int\nolimits_{t_p}^{t_{p+1}}a_{ij}(s)\,ds<\infty.
\ee
\end{defn}

Notice that in the discrete-time case the diagonal entries $a_{ii}(t)$ have no effect on the supremum $\ell$.

\begin{defn}\label{def.aqsc}
The matrix function $A(\cdot)$ is \emph{aperiodically} quasi-strongly connected (AQSC) if there exist $\ve>0$ and increasing sequence $t_p\to\infty$ satisfying the two conditions:
\begin{enumerate}[(i)]
\item the graph $\g[A_{t_p}^{t_{p+1}}]$ (defined in~\eqref{eq.unions}) is quasi-strongly $\ve$-connected for all $p$;
\item the function $A(\cdot)$ is $(t_p)$-bounded.
\end{enumerate}
\end{defn}

\begin{rem}
The term ``aperiodically'' emphasizes that the sequence $t_p$ can be \emph{arbitrary} unlike the usual uniform quasi-strong connectivity (where $t_{p+1}-t_p=T>0$). Durations $t_{p+1}-t_p>0$ can be arbitrarily large, so $(t_p)$-boundedness neither implies Assumption~\ref{asm.int-bound} nor follows from it.
\end{rem}
\begin{rem}\label{rem.ell-infty}
The supremum $\ell$ in~\eqref{eq.ell} is finite yet can be arbitrarily large. As $\ell\to\infty$, the AQSC condition relaxes to the quasi-strong connectivity of $\g_{\infty}$ (see Lemma~\ref{prop.necess}), which
condition is \emph{necessary} for consensus even when $\bar h=0$.
\end{rem}

\subsection{Consensus criteria}

We are now ready to formulate our first sufficient condition for consensus. We start with the continuous-time case.
\begin{thm}\label{thm.repeated-cont}
Suppose that algorithm~\eqref{eq.conse1d} satisfies Assumption~\ref{asm.int-bound} and that the matrix-valued function $A(\cdot)$ is AQSC.
Then, the algorithm establishes global consensus, no matter how large the delay bound $\bar h>0$ is.
Furthermore, a number $\theta\in (0,1)$ exists (depending only on $A(\cdot)$ and $h_{ij}(\cdot)$)
such that
\be\label{eq.contract}
\begin{gathered}
\|x(t)-\bar x\|_{\infty}\leq\La(t)-\la(t)\leq\theta^{k}(\La(t_*)-\la(t_*)),\\
\forall k\geq 0\quad\forall t:\, t_{r+2k(n-1)}\leq t\leq t_{r+2(k+1)(n-1)}.
\end{gathered}
\ee
Here $r$ is  an index such that $t_r\geq t_*$ and $n=|\V|$.
When $\bar h=0$ (undelayed case), Assumption~\ref{asm.int-bound} can be discarded.
\end{thm}

The conditions for the discrete-time case are similar, however, Assumption~\ref{asm.int-bound} is not needed and $A(t)$ is a stochastic matrix satisfying the \emph{strong aperiodicity} condition, which is typical for consensus criteria~\cite{ShiJohansson:13,XiaShiCao2019,ProCalaCao:2020}.
\begin{defn}
The sequence of stochastic matrices $A(t)$ is \emph{strongly aperiodic} if $a_{ii}(t)\geq\eta>0\,\forall i\in\V$.
\end{defn}
\begin{thm}\label{thm.repeated-disc}
Let the sequence of stochastic matrices $A(\cdot)$ be AQSC and strongly aperiodic. Then, algorithm~\eqref{eq.conse0d} establishes consensus, no matter how large the delay bound $\bar h>0$ is, and
~\eqref{eq.contract} retains its validity.
\end{thm}

As shown by Remark~\ref{rem.ell-infty}, the AQSC condition, being sufficient for the global consensus, is in fact ``very close'' to necessity: the necessary condition appears to be, in some sense, the ``limit case'' of
AQSC as $\ell\to\infty$ in~\eqref{eq.ell}.

In the case of a static matrix $A(t)\equiv A$, the strong aperiodicity condition guarantees that the graph $\g[A]$ is aperiodic (the inverse statement, however, is not true).
In the case when both weights $a_{ij}$ and delays are constant, the latter condition sometimes can be relaxed; e.g., consensus is established when $\g[A]$ is strongly connected and
contains at least one self-loop~\cite{ProCalaCao:2020}.  In the undelayed case, the strong periodicity condition can be replaced by other aperiodicity conditions coming from the inhomogeneous Markov chain theory~\cite{Bolouki:2016} whose verification, however, is a self-standing non-trivial problem.

Notice that the condition of $(t_p)$-boundedness in Definition~\ref{def.aqsc} is essential and cannot be discarded even if the graphs $\g[A_{t_1}^{t_2}]$ are complete and delays are absent $\bar h=0$, as shown by the counterexample in~\cite[Section~IV-C]{Moro:05}; the latter example deals with the discrete-time system, but its extension to the continuous time is straightforward.

\subsection{Alternative consensus conditions}

In this subsection, we discuss the relations between Theorems~\ref{thm.repeated-cont} and~\ref{thm.repeated-disc} and the previously published consensus criteria.

\subsubsection{AQSC vs. uniform quasi-strong connectivity}

As has been already mentioned, the \emph{uniform} quasi-strong connectivity (UQSC)~\cite{Moro:04,Sepul:09} is a special case of AQSC, where $t_p=pT$; in this case $(t_p)$-boundedness, obviously, is equivalent to Assumption~\ref{asm.int-bound},
and~\eqref{eq.contract} ensures \emph{exponentially fast} convergence to the consensus manifold $\{x: x_i=x_j\,\forall i,j\in\V\}$. Notice that in the literature some stronger versions of the UQSC property can be found, which require the matrix $A(t)$ to be uniformly bounded and piecewise constant with positive dwell time between subsequent switches~\cite{RenBeardBook,Muenz:11}. Theorem~\ref{thm.repeated-cont} discards the latter restrictions. In the discrete-time case, the coefficients $a_{ij}(t)$ are often supposed to be either $0$ or uniformly positive~\cite{Blondel:05,RenBeardBook}; we also discard this restriction for the off-diagonal entries $a_{ij}(t)$.
Theorems~\ref{thm.repeated-cont} and~\ref{thm.repeated-disc} can also be generalized also to \emph{nonlinear}
algorithms in the multidimensional vector spaces, as will be shown in Section~\ref{sec.ext}.

\subsubsection{Intermittent communication}

The difference between AQSC and UQSC properties is prominently illustrated by networks with \emph{intermittent} communication~\cite{Guanghui:12}:
\be\label{eq.intermittent}
a_{ij}(t)=\alpha(t)\bar a_{ij},\quad\forall i\ne j\quad\forall t\geq 0.
\ee
where $\bar a_{ij}\geq 0$ are \emph{constant} and $\alpha(t)\in\{0,1\}$ is an arbitrary measurable function, indicating the availability of the communication network at time $t$. Obviously, the persistent graph $\g_{\infty}$ either has no arcs (when $\alpha\in L_1[0,\infty)$) or coincides with $\g[\bar A]$ when $\int_0^{\infty}\alpha(t)dt=\infty$ (that is, the total time of network's availability is infinite).
\begin{lem}\label{prop.intermittent}
Assume that the graph $\g[\bar A]$ is quasi-strongly connected. The matrix function $A(\cdot)$ defined~\eqref{eq.intermittent} is AQSC if and only if $\int_0^{\infty}\alpha(t)dt=\infty$. Hence, the necessary consensus condition from Lemma~\ref{lem.necess}, in this situation, is also \emph{sufficient} for consensus under arbitrary bounded delays.
\end{lem}
\begin{IEEEproof}
Choosing arbitrary $\vk>0$, let $t_0=0$ and $t_{p+1}>t_p$ is defined in such a way that $\int_{t_p}^{t_{p+1}}\alpha(s)\,ds=\vk$. Obviously, the matrix $A(\cdot)$ is $(t_p)$-bounded and satisfies thus the AQSC condition (with $\ve=\vk\min\{\bar a_{ij}:\bar a_{ij}>0\}$). This proves the ``if'' part; the ``only if'' part is obvious.
\end{IEEEproof}

Notice that the condition from Lemma~\ref{prop.intermittent} allows arbitrarily long periods of ``silence'' when the network is unavailable. Unlike it, the UQSC condition requires that the network is periodically available: formally, $\int_t^{t+T}\alpha(s)\,ds\geq\ve>0$ for some $\ve>0$. As discussed in~\cite{Arcak:07,BarabanovOrtega2018}, the latter condition is similar in spirit to \emph{periodic excitation} in multi-agent control.

\subsubsection{The ``arc-balance'' condition}

The network from Lemma~\ref{prop.intermittent} is actually a particular case of a so-called \emph{arc-balanced} network~\cite{ShiJohansson:13-1,XiaShiCao2019}.
For simplicity, we consider only the ``anytime'' arc-balance from~\cite{ShiJohansson:13-1}; the Lemma below remains valid for the ``non-instantaneous'' arc-balance~\cite{XiaShiCao2019}.
\begin{defn}
The weighted graph $\g[A(\cdot)]$ is said to be \emph{arc-balanced} if a constant $K\geq 1$ such that
\be\label{eq.arc-balance}
K^{-1}a_{km}(t)\leq a_{ij}(t)\leq Ka_{km}(t)
\ee
for each pair of \emph{persistent} arcs $(m,k),\,(j,i)\in\E_{\infty}$.
\end{defn}

The main disadvantage of the arc balance condition is the necessity to know the persistent graph in advance. If $\g_{\infty}$ is quasi-strongly connected, then the arc balance entails
the AQSC, as shown by the following.
\begin{lem}\label{prop.arc-balance}
If persistent graph $\g_{\infty}$ is quasi-strongly connected, then the arc-balance condition entails the AQSC.
\end{lem}
\begin{IEEEproof}
Consider the continuous-time case, the discrete-time is analogous. Choose an arc $(j,i)\in\g_{\infty}$
and $\ve>0$. Since $a_{ij}\not\in L_1$, a sequence $t_p\to\infty$ exists such that
\[
\int_{t_p}^{t_{p+1}}a_{ij}(t)dt=K\ve.
\]
Then, using~\eqref{eq.arc-balance} we have
\[
\ve\leq\int_{t_p}^{t_{p+1}}a_{km}(t)dt\leq K^2\ve\quad\forall (m,k)\in\E_{\infty}.
\]
Recalling that $a_{k\ell}\in L_1$ for $(k,\ell)\not\in\E_{\infty}$, one proves $(t_p)$-boundedness of the matrix $A(\cdot)$. The graphs $\g[A_{t_p}^{t_{p+1}}]$ are quasi-strongly $\ve$-connected (their $\ve$-skeletons contain $\g_{\infty}$).
\end{IEEEproof}

Theorems~\ref{thm.repeated-cont} and~\ref{thm.repeated-disc} thus extend the consensus criteria from~\cite{ShiJohansson:13-1} (Theorems~3.1 and~4.1) to the case of delayed communication; as was already mentioned, they can be extended to the case of non-instantaneous arc-balance from~\cite{XiaShiCao2019}.

\section{The second sufficient condition: type-symmetry}\label{sec.symm}

Although, as shown by Remark~\ref{rem.ell-infty}, the AQSC conditions is rather close to necessity, some gap still remains between necessary and sufficient conditions. In the undelayed case, this gap can be removed provided that the interactions between the agents are \emph{reciprocal}, e.g. the weights are symmetric~\cite{CaoZheng:11}, type-symmetric~\cite{Lorenz:2005,Blondel:05} or cut-balanced~\cite{TsiTsi:13,Bolouki2015}; the most general conditions seem to be the non-instantaneous cut-balance~\cite{MartinHendrickx:2016,XiaShiCao2019} and weak reciprocity~\cite{ProCalaCao:2020}.
Most of the aforementioned results, in fact, are confined to the undelayed communication case; the only exceptions, to the best authors' knowledge, are the results for discrete-time consensus algorithms~\eqref{eq.conse0d} established in~\cite{Blondel:05} and~\cite{ProCalaCao:2020}; in both situations, however, a restrictive assumption is imposed on the weights $a_{ij}$, that is, the positive weights should be uniformly positive.

\begin{defn}\label{def.symm0}
 Matrix function $A(\cdot):[0,\infty)\to\r^{\V\times\V}$ with entries $a_{ij}\geq 0$ is \emph{type-symmetric} if $K\geq 1$ exists such that
\be\label{eq.type-symm}
K^{-1}a_{ji}(t)\leq a_{ij}(t)\leq Ka_{ji}(t)\quad\forall i,j\in \V\,\forall t\geq 0.
\ee
\end{defn}

Obviously, symmetric matrix ($A(t)=A(t)^{\top}\,\forall t\geq 0$) is type-symmetric (with $K=1$). A generalization of~\eqref{eq.type-symm} is the \emph{non-instantaneous}  type-symmetry introduced in~\cite{MartinHendrickx:2016}.
\begin{defn}\label{def.symm}
The matrix function $A(\cdot):[0,\infty)\to\r^{\V\times\V}$ with nonnegative entries $a_{ij}\geq 0$ possess the \emph{non-instantaneous type-symmetry} (NITS) property if
there exist an \emph{increasing} sequence $t_p\xrightarrow[p\to\infty]{}\infty$ (where $p=0,1,\ldots$) and a constant $K\geq 1$ such that
the following two conditions hold:
\begin{enumerate}
\item for any $i,j=1,\ldots,n$ one has
\begin{equation}\label{eq.type-symm-non}
\frac{1}{K}\int\limits_{t_p}^{t_{p+1}}a_{ji}(t)\,dt\leq \int\limits_{t_p}^{t_{p+1}}a_{ij}(t)\,dt\leq K\int\limits_{t_p}^{t_{p+1}}a_{ji}(t)\,dt;
\end{equation}
\item the function $A(\cdot)$ is $(t_p)$-bounded.
\end{enumerate}
 The definition for the discrete-time matrix function $A(t),\,t\geq 0$ remains same, replacing the integral in~\eqref{eq.type-symm-non} by the sum (in this case, $t_p$ are required to be integer).
\end{defn}
\begin{rem}
The condition~\eqref{eq.type-symm-non} implies, obviously, that the persistent graph $\g_{\infty}$ is undirected ($(i,j)\in\E_{\infty}$ if and only if $(j,i)\in\E_{\infty}$). In particular, it is quasi-strongly connected
if and only if it is strongly connected; otherwise, the graph $\g_{\infty}$ consists of several connected components.
\end{rem}

Verification of~\eqref{eq.type-symm-non} may seem a non-trivial problem, because the sequence $t_p$ should be same for all pairs $(i,j)$. In reality, this condition can be efficiently tested if the weights $a_{ij}$ are uniformly bounded (which also entails Assumption~\ref{asm.int-bound}), as follows from the proof of Theorem~2 in~\cite{MartinHendrickx:2016} (see also a practical example from~\cite{MartinHendrickx:2016}).

We now formulate criteria for consensus under the NITS condition in the continuous-time and discrete-time cases.
\begin{thm}\label{thm.symm-c}
Assume that the matrix-valued function $A(\cdot)$  with nonnegative entries $a_{ij}(t)\geq 0$ possess the NITS property and satisfies Assumption~\ref{asm.int-bound}. Then, every solution to~\eqref{eq.conse1d} converges: the limits~\eqref{eq.barx-i} exist for all $i\in\V$. The algorithm establishes consensus between agents $i$ and $j$ if and only if $i$ and $j$ belong to the same connected component of $\g_{\infty}$
(equivalently, a walk between $i$ and $j$ exists). The algorithm establishes global consensus \emph{if and only if} $\g_{\infty}$ is connected.
\end{thm}
\begin{thm}\label{thm.symm-d}
Assume that the sequence $\{A(t)\}$ of stochastic matrices possess the NITS property and is strongly aperiodic. Then, every solution to~\eqref{eq.conse0d} converges, and the algorithm establishes consensus
between agents $i$ and $j$ if they belong to the same connected component of $\g_{\infty}$.
\end{thm}

As illustrated by Proposition~3 in~\cite{MartinHendrickx:2016}, even in the undelayed case  ($\bar h=0$)  the condition of $(t_p)$-boundedness is essential and cannot be discarded (yet can be
replaced by some alternative restrictions on the matrix $A(\cdot)$~\cite{MartinGirard:2013,BarabanovOrtega2018}). However, in the special case of type-symmetric~\eqref{eq.type-symm} matrix,
this condition in fact is unnecessary as shown by the following.
\begin{rem}
In the discrete-time case ($A(t)$ is stochastic),~\eqref{eq.type-symm} entails the NITS condition. In the continuous-time case, the NITS condition follows from~\eqref{eq.type-symm} and Assumption~\ref{asm.int-bound}.
Indeed, choosing $t_p=p$, the condition~\eqref{eq.type-symm-non} follows from~\eqref{eq.type-symm}. The $(t_p)$-boundedness follows either from stochasticity (discrete time) or Assumption~\ref{asm.int-bound} (continuous time).
\end{rem}

\subsection{Comparison to alternative criteria}

The result of Theorem~\ref{thm.symm-c} extends the results of~\cite{CaoZheng:11,MatvPro:2013} on convergence of type-symmetric continuous-time algorithms to the delayed case. Notice however that, unlike those results, it does not guarantee  that $\dot x_i\in L_1[0,\infty)\,\forall i\in\V$; the validity of the latter statement remains in fact an open problem.
Theorem~\ref{thm.symm-d} generalizes the results on convergence of type-symmetric discrete-time algorithms from~\cite{Lorenz:2005,Blondel:05,Moro:05}; it should be noted that the result from~\cite{Blondel:05} (based on the earlier criteria from~\cite{Tsitsiklis}), in fact, allows bounded delays, however, assumes that nonzero weights $a_{ij}\ne 0$ are uniformly positive; the latter assumption is discarded in Theorem~\ref{thm.symm-d}.

In the undelayed case ($\bar h=0$), the type-symmetry condition can be further relaxed to \emph{cut-balance} and its non-instantaneous versions~\cite{TsiTsi:13,Bolouki2015,MartinHendrickx:2016,XiaShiCao2019}.
The relevant extension to the delayed case, however, remains a non-trivial problem.

Theorem~\ref{thm.symm-c} substantially differs from other results on delay-robust consensus: as shown\footnote{In~\cite{ProCala_CDC:2020}, the weights $a_{ij}$ are supposed to be bounded, however, the proof can be
generalized to matrices obeying Assumption~\ref{asm.int-bound}} in~\cite{ProCala_CDC:2020}, it retains its validity for the \emph{differential inequalities}
\be\label{eq.conse1d-ineq}
\dot x_i(t)\leq\sum\nolimits_{j\in\V}a_{ij}(t)(\hat x_j^i(t)-x_i(t)),\,\forall i\in\V,\,t\geq 0.
\ee
Whereas the consensus dynamics  has some well-known contraction properties (the diameter $\La(t)-\la(t)$ of the convex hull,
spanned by the agents' values can serve as a Lyapunov function), inequality~\eqref{eq.conse1d-ineq}
does not possess such a properties: whereas $\La(t)$ is non-increasing, $\la(t)$ need not be monotone.
In particular, Theorem~\ref{thm.repeated-cont} is not valid for the inequalities and estimates like~\eqref{eq.contract}
cannot be established for them (in fact,
consensus in~\eqref{eq.conse1d-ineq} requires \emph{strong} connectivity of the persistent
graph $\g_{\infty}$~\cite{ProCao:2017,ProCala_CDC:2020}).
Theorem~\ref{thm.symm-c} thus requires some tools that are principally different from usual contraction analysis~\cite{Bliman:06,Muenz:11}; its proof is actually based on
the seminal idea of the solution's ordering originally used to study undelayed type-symmetric algorithms in~\cite{MatvPro:2013} (a similar yet slightly different technique was
employed in~\cite{TsiTsi:13,Bolouki2015}). As shown in our previous works~\cite{ProCao:2017,ProCalaCao:2020}, averaging inequalities, which are beyond the scope of this paper, have numerous applications in
multi-agent control and social dynamics modeling.

%%%%%%%%5
\section{Further extensions}\label{sec.ext}

In this section, we discuss some extensions and applications of the main results. We show that consensus criteria remain
valid in presence of vanishing exogenous disturbances (Subsect.~\ref{subsec.robust}) and can be generalized, without
principal differences, to nonlinear averaging algorithms (Subsect.~\ref{subsec.nlin}). We also show that
the first pair of consensus criteria (Theorems~\ref{thm.repeated-cont} and~\ref{thm.repeated-disc}) can be applied to analyze
algorithms known as \emph{containment control} with multiple leaders and \emph{target aggregation}.

\subsection{Robustness against disturbances}\label{subsec.robust}

In this section, we consider robustness of consensus against exogenous \emph{disturbances}, examining the systems\footnote{The disturbed discrete-time consensus algorithms have same
properties as the system~\eqref{eq.conse1d-f} (and can be reduced to~\eqref{eq.conse1d-f}  using the trick exploited
in Subsect.~\ref{subsec.proof-discrete}), we thus do not consider them for the sake of brevity.}
\begin{gather}
\dot x_i=\sum\limits_{j\in\V}a_{ij}(t)(\hat x_j^i(t)-x_i(t))+f_i(t),\,i\in\V\label{eq.conse1d-f}
\end{gather}
Functions $f_i$ are locally $L_1$-summable and the solution is determined by specifying initial condition~\eqref{eq.initial}.

The first lemma, which is very general and does not require even Assumption~\ref{asm.int-bound}, states that consensus is always robust to $L_1$-summable disturbances.
\begin{lem}\label{lem-l1-rob}
If the undisturbed algorithm~\eqref{eq.conse1d} establishes consensus between two agents $i$ and $j$, the same
holds for~\eqref{eq.conse1d-f} with an $L_1$-summable disturbance: $\int_0^{\infty}|f_m(t)|<\infty\,\forall m$.
\end{lem}

If the AQSC condition holds, then consensus is robust against a much broader class of disturbances.
\begin{lem}\label{lem-rob-general}
Suppose that the assumptions of Theorem~\ref{thm.repeated-cont} hold. Then, for any solution to~\eqref{eq.conse1d-f}, the inequality
holds as follows
\be\label{eq.rob-general}
\varlimsup_{t\to\infty}(\max x(t)-\min x(t))\leq C\varlimsup_{p\to\infty}\int\limits_{t_p}^{t_{p+1}}\|f(t)\|_{\infty}dt,
\ee
where $C=C(\theta)$ is a constant determined by $\theta$ from~\eqref{eq.contract}.
In particular, algorithm~\eqref{eq.conse1d-f} establishes global consensus when
\be\label{eq.vanishing}
\int_{t_p}^{t_{p+1}}|f_i(t)|\,dt\xrightarrow[p\to\infty]{} 0\quad\forall i\in\V.
\ee
\end{lem}
\begin{rem}
In the special case where the stronger \emph{uniform} quasi-strong connectivity holds ($t_p=pT$), condition~\eqref{eq.vanishing}
holds e.g. when the disturbance is \emph{vanishing} $f(t)\xrightarrow[t\to\infty]{}0$. The latter result in the undelayed
case ($\bar h=0$) was first established in~\cite{ShiJohansson:13}; it was also shown that the UQSC condition is
\emph{necessary} for this type of robustness.
\end{rem}

%\subsection{Consensus between identical neutrally stable agents}

\subsection{Nonlinear multidimensional consensus dynamics}\label{subsec.nlin}

Up to now, we have considered agents with scalar values $x_i\in\r$. The results trivially extend to the multidimensional case
where $x_i\in\r^m$, because algorithms~\eqref{eq.conse1d} and~\eqref{eq.conse0d} alter each dimension independently.
An important property of such algorithms is the non-expansion of convex
hull~\cite{Muenz:11}, which can be formulated as follows.
\begin{lem}\label{lem.convex}
Consider algorithm~\eqref{eq.conse0d} or~\eqref{eq.conse1d} with $x_i\in\r^m$.
Let $\mathcal{D}\subseteq\r^m$ be a convex set and $x_i(t)\in\mathcal{D}\,\forall i\in\V$ when $t_*-\bar h\leq t\leq t_*$.
Then, $x_i(t)\in\mathcal{D}$ for all $t\geq t_*$.
\end{lem}

Furthermore, the results can be extended to a broad class of \emph{nonlinear} consensus algorithms similar to those
considered in~\cite{Bliman:06,Antonis,Muenz:11}.

In particular, we may consider the counterpart of~\eqref{eq.conse1d}
\be\label{eq.conse1nd}
\dot x_i(t)=\sum_{j\ne i}a_{ij}(t)\vp_{ij}(\hat x^i_j(t),x_i(t))\in\r^m,
\ee
where nonlinear coupling functions $\vp_{ij}$ are such that
\be\label{eq.phi-struct}
\vp_{ij}(y,z)=\psi_{ij}(y,z)(y-z),\;\; 0<\psi_{ij}(y,z)\in\r\,\forall y,z\in\r^m.
\ee

\begin{thm}\label{thm.nlin-c}
Theorems~\ref{thm.repeated-cont} and~\ref{thm.symm-c} retain their validity for algorithm~\eqref{eq.conse1nd}, provided
that $\vp_{ij}$ are continuous, structured as in~\eqref{eq.phi-struct}, and for each compact set $\Omega\subset\r^{2m}$
one has
\be\label{eq.psi-ineq}
0<\inf_{(y,z)\in\Omega}\psi_{ij}(y,z)\leq \sup_{(y,z)\in\Omega}\psi_{ij}(y,z)<\infty
\ee
(this holds e.g. when $\psi_{ij}>0$ are continuous functions).
\end{thm}
\vskip0.1mm

Theorem~\ref{thm.nlin-c} generalizes the result of~\cite[Theorem~1]{Muenz:11} in several directions: first, the UQSC assumption
is relaxed (to either AQSC or the NITS with integral connectivity), second, the structure of nonlinearities is more general
than in~\cite{Muenz:11}, where it was assumed
that $\underline\varkappa\|y-z\|\leq\psi_{ij}(y,z)\leq\overline\varkappa\|y-z\|$ for some constants
$\underline\varkappa,\overline\varkappa>0$.

Similarly properties hold for the discrete-time dynamics
\be\label{eq.conse0nd}
x_i(t+1)=x_i(t)+\sum_{j\ne i}a_{ij}(t)\vp_{ij}(\hat x^i_j(t),x_i(t))\in\r^m,
\ee
where $\phi_{ij}$ have structure~\eqref{eq.phi-struct} with $\psi_{ij}\in(0,1]$.

\begin{thm}\label{thm.nlin-d}
Theorems~\ref{thm.repeated-disc} and~\ref{thm.symm-d} retain their validity for algorithm~\eqref{eq.conse0nd}, provided
that $\vp_{ij}$ are structured as in~\eqref{eq.phi-struct}, $\psi_{ij}\in(0,1]$ and for each
compact set $\Omega\subset\r^{2m}$ the leftmost inequality in~\eqref{eq.psi-ineq} holds.
\end{thm}
\vskip0.1mm

\subsection{Containment control and target aggregation}\label{subsec.contain}

In this subsection, we consider some applications to control of mobile agents that are known as
\emph{containment control}~\cite{RenCaoBook} and \emph{target aggregation}~\cite{ShiHong:09}.
Usually, these algorithms are considered in continuous-time, although the discrete-time counterpart of
results below can also be obtained.

We start, however, with a formally simpler dynamics known as the ``leader-following'' consensus algorithm.

\subsubsection{Damped averaging dynamics and leader-following consensus}

Consider the following modification of~\eqref{eq.conse1nd}
\begin{gather}
\dot x_i(t)=-d_i(t)x_i(t)+\sum_{j\in\V}a_{ij}(t)(\hat x_j^i(t)-x_i(t))\label{eq.avg-cont0-damp}.
\end{gather}
As usual, matrix $A(t)$ satisfies Assumptions~\ref{asm.coeff} and $d_i(t)\geq 0$.

Usually, equations~~\eqref{eq.avg-cont0-damp} arise in regard to a formally more general algorithm of consensus
with static\footnote{In the case of time-varying ``leader'', the consensus problem considered below
is known as the reference-tracking consensus~\cite{RenBeardBook}. This problem is beyond the scope of this paper.
If the trajectory of the leader and its derivative are unknown, one may guarantee approximate
consensus (using e.g. robustness results from Section~\ref{subsec.robust}) or use nonlinear
 sliding mode control~\cite{LewisBook} to guarantee asymptotic convergence to the desired trajectory.}
 ``leader''~\cite{RenBeardBook} whose index is $\omega\not\in\V$ and whose value $x_{\omega}$ is constant.
  At each time instant, only some agents are able to access the leader's value , such agents can be called ``informed''.
 In order to make all agents reach consensus with the leader, the following generalization of system~\eqref{eq.avg-cont0-damp} may be considered:
\begin{gather}
\dot x_i(t)=d_i(t)(x_{\omega}-x_i(t))+\sum_{j\in\V}a_{ij}(t)(\hat x_j^i(t)-x_i(t))\label{eq.avg-cont0-lead}.
\end{gather}
The informed (at time $t$) agents correspond to $d_i(t)>0$.

The following proposition is straightforward, using the transformation $x_i\mapsto x_i-x_{\omega}$ and noticing that $x_{\omega}={\rm const}$.
\begin{prop}\label{prop.leader-stab}
Algorithm~\eqref{eq.avg-cont0-lead} establishes (global) \emph{leader-following consensus}
\be\label{eq.leader-foll-conse}
\lim_{t\to\infty}x_i(t)=x_{\omega}\quad\forall i\in\V
\ee
if and only if the system~\eqref{eq.avg-cont0-damp} is asymptotically stable.
\end{prop}

Considering the augmented matrix $\hat A$, where
\be\label{eq.abar}
\hat a_{ij}=
\begin{cases}
a_{ij}, &i,j\in\V,\\
d_i, &i\in\V,j=\omega,\\
0, &i=\omega,j\in\V,\\
0, &i=j=\omega,
\end{cases}
\ee
system~\eqref{eq.avg-cont0-lead} becomes a special case of averaging
algorithm~\eqref{eq.conse1d}.

In the graph $\g[\hat A(\cdot)]$, the leader's node $\omega$ has no incoming arcs (but for the self-loop), that is,
the leader is independent of the remaining agents yet influences them. Due to this influence, the graph can
be AQSC, whereas the NITS condition can never be guaranteed.
From Lemma~\ref{lem.necess} and Theorem~\ref{thm.repeated-cont}, the following corollary is immediate.
\begin{cor}\label{cor.lead-foll-repeated}
For asymptotic stability of system~~\eqref{eq.avg-cont0-damp} it is necessary that the persistent graph $\hat\g_{\infty}$
of the \textbf{augmented} matrix function $\hat A(\cdot)$
is quasi-strongly connected. System~\eqref{eq.avg-cont0-damp} is stable if $\hat A(\cdot)$ satisfies the conditions
of Theorem~\ref{thm.repeated-cont}.
\end{cor}

Notice that the AQSC condition requires that implies that the leader's influence on the other agents is ``non-negligible''
over each interval $[t_p,t_{p+1}]$ in the sense that
\[
\sum_{i\in\V}\int_{t_p}^{t_{p+1}}d_i(s)\,ds\geq\ve>0.
\]

\subsubsection{Containment control with multiple leaders and target aggregation}

A natural extension of the previously considered algorithm~\eqref{eq.avg-cont0-lead} is the algorithm of \emph{containment control} with multiple ``leaders''.
In this problems, the agents are usually considered as mobile objects (robots, vehicles) moving in some finite-dimensional
space $\mathbb{X}=\r^m$; the values $x_i\in\mathbb{X}$ stand for positions of the objects (and are also referred to as
state vectors). The leaders are dedicated agents indexed by the set $\V_{L}$ ($L$=leader) whose dynamics can be arbitrary
and do not depend on the remaining agents (called ``followers''). In this subsection, we confine ourselves to containment
control with \emph{static} leaders; different approaches to containment control with dynamics leaders can be found
in~\cite{RenCaoBook,LewisBook}. The goal of a containment control algorithm is to bring all followers into the convex
hull spanned by the leaders' positions. It is remarkable that, long before control theorists, the algorithms of
containment control were proposed by mathematical sociologists~\cite{ProTempo:2017-1}.

The most standard algorithm of containment control mimics the averaging dynamics: the followers move according to
\begin{gather}
\dot x_i(t)=\sum_{j\in\V}a_{ij}(t)(\hat x_j^i(t)-x_i(t))+\sum_{k\in\V_{L}}b_{ik}(t)(x_{k}-x_i(t)).\label{eq.avg-cont0-contain}
\end{gather}
Here $i\in\V$, $x_i(t),x_{k}\in\mathbb{X}$, $B\in\r^{\V\times\V_L}$ stands for a nonnegative matrix, $b_{ik}$ measures ``attraction'' of follower $i$ to leader $\ell$. In the continuous-time case, $A(t)\in\r^{\V\times\V}$ is an arbitrary nonnegative matrix.
We denote
\be\label{eq.aux-d-i}
d_i(t)\dfb \sum_{k\in\V_L}b_{ik}(t)\quad\forall i\in\V.
\ee

%Similar to the leader-following consensus algorithms from the previous subsection, the containment algorithms
%may be considered as a special class of averaging dynamics~\eqref{eq.conse1d} by introducing the augmented
%matrix $\tilde A=(\tilde a_{ij})$
%\be\label{eq.atilde}
%\tilde a_{ij}=
%\begin{cases}
%a_{ij}, &i,j\in\V,\\
%b_{ij}, &i\in\V,j\in\V_L,\\
%0, &i\in\V_L,j\in\V,\\
%1, &i,j\in\V_L.
%\end{cases}
%\ee

\begin{defn}
Algorithm~~\eqref{eq.avg-cont0-contain} provides \emph{containment} if all followers converge
asymptotically to the convex hull spanned by the leaders' states,\footnote{Here $\dist$ stands for the Euclidean distance from
a point to a set.} that is,
\be\label{eq.containment}
\dist\{x_i(t),\mathcal{C}_L\}\xrightarrow[t\to\infty]{}0\,\forall i\in\V,\; \mathcal{C}_L\dfb\conv\{x_{k}\}_{k\in\V_L}
\ee
for any initial condition and any positions of the leaders.
\end{defn}

It appears that in fact the containment property is \emph{equivalent} to the stability of damped averaging dynamics,
studied in the previous subsection. Hence, convergence of containment control algorithm, in fact, does not require special
theory and can be derived from conventional consensus criteria (Corollary~\ref{cor.lead-foll-repeated}).
\begin{lem}\label{Lem.contain-stab}
Algorithm~\eqref{eq.avg-cont0-contain} provides containment if and only if
damped averaging dynamics~\eqref{eq.avg-cont0-damp}, where $d_i(t)$ is defined in~\eqref{eq.aux-d-i},
  are asymptotically stable.
  %\item algorithm~\eqref{eq.avg-cont0-lead} establishes leader-follower consensus.
%\end{enumerate}
\end{lem}

Combining Corollary~\ref{cor.lead-foll-repeated} and Lemma~\ref{Lem.contain-stab}, the following result is immediate.
\begin{cor}\label{cor.lead-foll-contain}
Algorithm~\eqref{eq.avg-cont0-contain} provides containment if matrix $\hat A(\cdot)$ defined in~\eqref{eq.abar} satisfies the conditions
of Theorem~\ref{thm.repeated-cont}.
\end{cor}

Comparing the result of Corollary~\ref{cor.lead-foll-repeated} to standard criteria of convergence,
e.g.~\cite[Lemma~5.6]{RenCaoBook}, one notices that a number of restrictions have been relaxed, in particular:
the matrices $A(t)$ and $B(t)$ need not be piecewise-continuous with positive dwell time between consecutive switchings;
there is no uniform positivity\footnote{Formally, Lemma~5.6 of~\cite{RenCaoBook} does not impose this assumption,
however, it is used implicitly in Lemma~1.5, on which the proof of Lemma~5.6 relies} of non-zero
entries: $a_{ij}\in\{0\}\cup[\eta,\infty)$; the assumption of UQSC\footnote{The uniform quasi-strong connectivity of $\hat A$ is equivalent to the property, called in~\cite{RenCaoBook} the existence of
``united directed spanning tree''.} for $\hat A$ is relaxed to the AQSC.

It should be noted that containment control algorithms often arise under different names, e.g. in the very recent
work~\cite{Alexandrov:19} a special algorithm with two static leaders ($|\V_L|=2$) was considered, which provides uniform deployment
of the follower agents on the line segment connecting the two leaders. Delay robustness of this algorithm can be derived
from Corollary~\ref{cor.lead-foll-repeated} (a direct proof given in~\cite{Alexandrov:19} assumes that delays are constant).

\subsubsection{Algorithms of target aggregation}

A natural extension of the containment control algorithm was proposed in~\cite{ShiHong:09} under name of \emph{target
aggregation}, or gathering of mobile agents in a convex closed set $\Omega\subset\mathbb{X}=\r^m$.
Unlike the containment control algorithm, this set need not be a convex polytope and the agents, in principle,
have no information about its structure. Some ``informed'' agents, however, at time $t$ may have access to
elements $\omega_i(t)\in\Omega$ of this set.

The target aggregation algorithm~\cite{ShiHong:09} resembles~\eqref{eq.avg-cont0-lead}
and has the following structure\footnote{In algorithm from~\cite{ShiHong:09}, $\omega_i(t)-x_i(t)$ can be
replaced by a more general nonlinearity $f_i(t,x_i)$, providing attraction to the set $\Omega$.
However, $f_i(t,x_i)=\omega_i(t)-x_i$ is the most typical class of functions, satisfying all the assumptions adopted
in~\cite{ShiHong:09}. Unlike~\cite{ShiHong:09}, our algorithm allows communication delays.}:
\begin{gather}
\dot x_i(t)=\sum_{j\in\V}a_{ij}(t)(\hat x_j^i(t)-x_i(t))+d_i(t)(\omega_i(t)-x_i(t))\in\mathbb{X}.\label{eq.avg-cont0-agreg}
\end{gather}
Here $i\in\V$, $\omega_i:[0,\infty)\to\Omega$ is Lebesgue measurable and locally $L_1$-summable in order to guarantee
the solution existence. Notice that if $d_i(t)=0$ (agent $i$ is not ``informed'' at time $t$), then $\omega_i(t)$
plays no role and can be arbitrary.

Similar to containment control, the goal of the algorithm is to gather all agents in the target set $\Omega$.
\begin{defn}
Algorithm~\eqref{eq.avg-cont0-contain} provides \emph{aggregation} of the agents in
target set $\Omega$ if for any initial condition their positions converge asymptotically to this set
\be\label{eq.aggreg}
\dist\{x_i(t),\Omega\}\xrightarrow[t\to\infty]{}0\quad\forall i\in\V.
\ee
\end{defn}

Similarly to Lemma~\ref{Lem.contain-stab}, one can prove the following.
\begin{lem}\label{Lem.aggreg-stab}
Algorithm~\eqref{eq.avg-cont0-agreg} provides target aggregation for \emph{any} choice of the set $\Omega$ and
functions $\omega_i(t)\in\Omega$ (abiding by the aforementioned conditions) if and only if
  the damped averaging dynamics~\eqref{eq.avg-cont0-damp} are asymptotically stable.
\end{lem}

Lemma~\ref{Lem.aggreg-stab} shows that linear target aggregation algorithms in fact also do not require to develop a special
theory, boiling down to the damped averaging dynamics.

\section{Technical proofs}\label{sec.proof}

The proofs rely on several technical lemmas that are collected in Subsect.~\ref{subsec.technical}.

%%%%%%%%%%%%
\subsection{Evolutionary matrices and their properties}\label{subsec.technical}

%Properties of the evolutionary matrices (proved in~\cite{ProCal:2020arx})}

Similar to linear ODE, a linear delay differential equation admits a so-called \emph{evolutionary matrix}~\cite{Hale,FridmanBook:2014} which will be denoted $U(t,t_*)$ (where $t\geq t_*\geq 0$).

By definition, the $i$th column of $U(t,t_*)$ is the solution of~\eqref{eq.conse1d} that corresponds to the initial condition $x(t_*)=\mathbf{e}_i$ (the $i$th basis vector) and $x(t)\equiv 0\,\forall t<t_0$
(recall that the solution is right-continuous at $t=t_*$). Obviously, the solution with $x(t_*)=x_*$ and $x(t)\equiv 0\,\forall t<t_0$ is then given by $U(t,t_*)x_*$.
For delay systems,  the Cauchy (``variation of constants'') formula retains its validity, stating that the solution to the \emph{disturbed} systems of equations~\eqref{eq.conse1d-f} obeys the equation
\be\label{eq.cauchy+}
\begin{gathered}
x(t)=U(t,t_*)x(t_*)+\int_{t_*}^tU(t,\xi)[f(\xi)+g(\xi)]\,d\xi,\\
g_i(t)\dfb \sum_{j\ne i}a_{ij}(t)x_j(t-h_{ij}(t))\chi_{t_*}(t-h_{ij}(t))\,\forall i\in\V,
\end{gathered}
\ee
where the indicator function $\chi_{t_*}$ is defined as
\be\label{eq.indicator}
\chi_{t_*}(t)\dfb\begin{cases}
0,\quad t\geq t_*,\\
1,\quad t<t_*.
\end{cases}
\ee
In particular, the solution with initial condition~\eqref{eq.initial} is found from~\eqref{eq.cauchy+} by substituting
$x(t_*)=x_*$ and $x_j(t-h_{ij}(t))=\vp(t-h_{ij}(t))$ for all $t$ such that $t-h_{ij}(t)<t_*$.

In this subsection, we collect a number of technical lemmas that are
concerned with the properties of evolutionary matrices of system~\eqref{eq.conse1d}.

It is well-known~\cite{Antonis} that the convex hulls spanned by the agents' values at time $t$ (in our situation, $[\la(t),\La(t)]$) are nested. We formulate a more general statement, which is used in analysis of
the evolutionary matrices.
 \begin{lem}\label{prop.bound}
The functions $\la(t),\La(t)$ are respectively  \emph{non-decreasing} and \emph{non-increasing} for any solution of~\eqref{eq.conse1d}. Furthermore, for any solution of~\eqref{eq.conse1d} the inequalities hold
\be\label{eq.simple-est}
\begin{gathered}
x_i(t)\leq x_i(s)e^{-\int_s^t\alpha_i(\xi)d\xi}+\La(s)\left[1-e^{-\int_s^t\alpha_i(\xi)d\xi}\right],\\
x_i(t)\geq x_i(s)e^{-\int_s^t\alpha_i(\xi)d\xi}+\la(s)\left[1-e^{-\int_s^t\alpha_i(\xi)d\xi}\right],\\
\alpha_i(t)\dfb\sum_{j\ne i} a_{ij}(t), \,\forall i\in\V\,\forall t\geq s\geq 0.
\end{gathered}
\ee
Finally, if Assumption~\ref{asm.int-bound} holds, then
\be\label{eq.lims}
\lim_{t\to\infty}\min x(t)=\lim_{t\to\infty}\la(t), \lim_{t\to\infty}\max x(t)=\lim_{t\to\infty}\La(t).
\ee
%The function $\La(t)$ is also non-increasing for any feasible solution of the inequality~\eqref{eq.conse1d-ineq} (e.g., all solutions of~\eqref{eq.conse1d}).
\end{lem}
\begin{IEEEproof}
We prove only the statements involving $\La(t)$, which remain valid for inequalities~\eqref{eq.conse1d-ineq}.
The statement related to $\la(t)$ (valid for equation~\eqref{eq.conse1d}) are trivially derived from them by considering the solution $(-x(t))$, which also obeys~\eqref{eq.conse1d}.

We first prove that $\La(t)$ is non-increasing for any feasible solution to inequality~\eqref{eq.conse1d-ineq}. For such a solution, choose arbitrary constant $\La'>\La(t_*)$. We are going to show that $\max x(t)<\La'\,\forall i\in\V$ for any $t\geq t_*$. Obviously, the latter inequality holds when $t$ is close to $t_*$; let $t'$ be the \emph{first} instant $t>t_*$ when it is violated, that is,
\[
\begin{gathered}
\max x(t)<\La'\quad\forall t\in [t_*-\bar h,t'),\\
x_{i}(t')=\La'\quad\text{for some $i$.}
\end{gathered}
\]
One arrives at a contradiction with~\eqref{eq.conse1d-ineq}, because
\be\label{eq.aux0}
\begin{gathered}
\dot x_i(t)\leq \alpha_i(t)[\La'-x_i(t)]\,\forall t\in [t_*,t')\Longrightarrow\\
x_i(t')\leq e^{-\int_{t_*}^{t'}\alpha_i(s)ds}x_i(t_*)+\La'\left(1-e^{-\int_{t_*}^{t'}\alpha_i(s)ds}\right)<\La'
\end{gathered}
\ee
(where $\alpha_i$ is defined in~\eqref{eq.simple-est}). The contradiction proves that $\La(t)<\La'\;\;\forall t\ge t_*$. Since $\La'>\La(t_*)$ can be arbitrary, we have $\La(t)\leq \La(t_*)$ whenever $t\geq t_*$. Replacing $t_*$ by an arbitrary $s\geq t_*$, one proves similarly that $\La(t)\leq\La(s)\,\forall t\geq s$.

The first inequality in~\eqref{eq.simple-est} is proved similar to~\eqref{eq.aux0}. Indeed, when $t\geq s$, we have $\hat x_j^i(t)\leq\La(t)\leq\La(s)\,\forall j\ne i$. and hence
\[
\dot x_i(t)\leq \alpha_i(t)[\La(s)-x_i(t)]\quad\forall t\geq s\geq t_*.
\]

If Assumption~\ref{asm.int-bound} holds, one has $\exp(-\int_{s}^t\alpha_i(\xi)d\xi)\geq\theta\dfb=e^{-(n-1)\mu_{\bar h}}$ whenever $t\in[s,s+\bar h]$, and thus
\[
\begin{gathered}
\La(s+\bar h)=\max\limits_{t\in[s,s+\bar h]}\max x(t)\leq \theta\max x(s)+(1-\theta)\La(s)\Longrightarrow\\
\theta^{-1}[\La(s+\bar h)-(1-\theta)\La(s)]\leq\max x(s)\leq\La(s),
\end{gathered}
\]
which proves the first inequality in~\eqref{eq.lims}.
\end{IEEEproof}

Lemma~\ref{prop.bound} implies the following important property of the evolutionary matrices.
 \begin{lem}\label{lem.cauchy0}
 Matrix $U(t,s)$ is \emph{substochastic}\footnote{In the case undelayed case ($\bar h=0$), the evolutionary matrix is known to be stochastic~\cite{RenBeardBook,Bolouki:2016}.
 This is the principal difference between the delayed and undelayed consensus dynamics, which makes it impossible to reduce delayed equations~\eqref{eq.conse1d} to ergodicity of Markov chains
 and stochastic matrix products.} for all $t\geq s\geq 0$.  For any solution of~\eqref{eq.conse1d-f}, one has
 \be\label{eq.cauchy-del-estim}
 \begin{gathered}
\la(t_*)[\ones-U(t,t_*)\ones]\leq \\x(t)-U(t,t_*)x(t_*)-\int_{t_*}^tU(t,\xi)f(\xi)d\xi\leq\\
\leq \La(t_*)[\ones-U(t,t_*)\ones].
\end{gathered}
\ee
 If Assumption~\ref{asm.int-bound} holds, then $U(t,s)$ have \emph{uniformly} positive row sums: a constant $\psi>0$ exists such that $U(t,s)\ones\geq\psi\ones$ whenever $t\geq s\geq 0$.
 \end{lem}
 \begin{IEEEproof}
 By construction of $U(t,t_*)$, the solution $x(t)=U(t,t_*)\mathbf{e}_i$ corresponds to $\la(t_*)=0$ and hence $x_i(t)\geq 0\,\forall i\in\V\,\forall t\geq t_*$. Therefore, $U(t,t_*)$ is a nonnegative matrix. Similarly, the solution $x(t)=U(t,t_*)\ones$ corresponds to $\La(t_*)=1$, and hence $U(t,t_*)\ones\leq\ones\,\forall t\geq t_*$, so $U(t,t_*)$ is substochastic whenever $t\geq t_*\geq 0$.

 Consider now the solution $x(t)\equiv 1\,\forall t\geq t_*$, which corresponds to the initial condition~\eqref{eq.initial} with $x_*=1$ and $\vp(s)=1\,\forall s\in[-\bar h,0)$. Using~\eqref{eq.cauchy+} with $f=0$, one has
 \[
 \begin{gathered}
    \ones=U(t,t_*)\ones+\int_{t_*}^tU(t,\xi)g^0(\xi)\,d\xi,\\
    g_i^0(t)\dfb \sum_{j\ne i}a_{ij}(t)\chi_{t_*}(t-h_{ij}(t))\,\forall i\in\V.
\end{gathered}
 \]
 For a general solution to~\eqref{eq.conse1d-f}, we have $\la(t_*)\leq x_i(t)\leq\La(t_*)$ whenever $t\in[t_*-\bar h,t_*]$, and thus the functions $g_i(t)$ introduced in~\eqref{eq.cauchy+} obey the inequalities $\la(t_*)g_i^0(t)\leq g_i(t)\leq \La(t_*)g_i^0(t)$. For this reason,
 \[
 \begin{gathered}
 x(t)-U(t,t_*)x(t_*)-\int_{t_*}^tU(t,\xi)f(\xi)d\xi\geq\\
 \geq \la(t_*)\int_{t_*}^tU(t,\xi)g^0(\xi)\,d\xi=\la(t_*)\left(\ones-U(t,t_*)\ones\right)
 \end{gathered}
 \]
 and, similarly,
 \[
 \begin{gathered}
 x(t)-U(t,t_*)x(t_*)-\int_{t_*}^tU(t,\xi)f(\xi)d\xi\leq\\
 \leq \La(t_*)\int_{t_*}^tU(t,\xi)g^0(\xi)\,d\xi=\La(t_*)\left(\ones-U(t,t_*)\ones\right),
 \end{gathered}
 \]
 which finishes the proof of~\eqref{eq.cauchy-del-estim}.

 To prove the final statement, the second inequality from~\eqref{eq.simple-est} can be applied to the solution $x(t)=U(t,s)\ones$. Obviously, for this solution one has $\la(s)=0$ and $x_i\geq 1\,\forall i\in\V$. Hence for any $t\in[s,s+\bar h]$ one has
 \[
 x_i(t)\overset{\eqref{eq.int-bound}}{\geq}\psi\dfb e^{-(n-1)\mu_{\bar h}}\;\forall i\in\V.
 \]
 In particular, $\la(s+\bar h)\geq\psi$, and thus (Lemma~\ref{prop.bound}) $x(t)\geq\psi\ones$ also for $t\geq s+\bar h$.
 This finishes the proof
 \end{IEEEproof}
  \begin{rem}\label{rem.monotone}
  Non-negativity of matrices $U(t,t_*)$ implies, in view of~\eqref{eq.cauchy+} in particular, the monotonicity of dynamics~\eqref{eq.conse1d}. If $x,\tilde x$ are two solutions such that $x(t)\leq\tilde\tilde x(t)$ for $t\in[t_*-\bar h,t_*]$, then $x(t)\leq\tilde x(t)$ for $t\geq t_*$.
 \end{rem}

Consensus, similar to the undelayed case~\cite{RenBeardBook}, admits a simple interpretation in terms of the evolutionary matrix.
\begin{lem}\label{lem.consensus-cauch}
The algorithm~\eqref{eq.conse1d} establishes consensus among agents $i$ and $j$ if and only if the limits exist and coincide:
\be\label{eq.conse-cauchy}
\lim_{t\to\infty}\mathbf{e}_i^{\top}U(t,t_*)=\lim_{t\to\infty}\mathbf{e}_j^{\top}U(t,t_*)\quad\forall t_*\geq 0.
\ee
Global consensus is established if and only if the limit exists
\be\label{eq.u-bar}
\bar U_{t_*}\dfb\lim_{t\to\infty}U(t,t_*)=\ones p^{\top}_{t_*},\quad p_{t_*}\in\r^{\V}
\ee
(that is, $\bar U_{t_*}$ has equal rows) for any $t_*\geq 0$.
\end{lem}
\begin{IEEEproof}
It suffices to prove the first statement, the second one follows from it.

The ``only if'' part in the first statement is straightforward, because $U(t,t_*)x_*$ is a solution to~\eqref{eq.conse1d} for each vector $x_*\in\mathbb{R}^{\V}$. Consensus between agents $i,j$ implies that $\lim_{t\to\infty}\lim_{t\to\infty}\mathbf{e}_i^{\top}U(t,t_*)x_*=\lim_{t\to\infty}\mathbf{e}_j^{\top}U(t,t_*)x_*$, which is equivalent to~\eqref{eq.conse-cauchy}.

To prove the ``if'' part, recall that an arbitrary solution to~\eqref{eq.conse1d} satisfies equation~\eqref{eq.cauchy+} with $f\equiv 0$. By noticing that $g(t)=0$ for $t>\bar h$, we have
\[
x(t)=U(t,t_*)x(t_*)+\int_{0}^{\bar h}U(\xi,t_*)g(\xi)dt\;\;\forall t>\bar h.
\]
Using the Lebesgue dominated convergence theorem, it can be easily shown that~\eqref{eq.conse-cauchy} entails
\[
\begin{aligned}
\bar x_i&=\lim_{t\to\infty}\mathbf{e}_i^{\top}x(t)=\lim_{t\to\infty}\mathbf{e}_i^{\top}U(t,t_*)x(t_*)+\\&+\int_{t_*}^{\bar h}\lim_{\xi\to\infty}\mathbf{e}_i^{\top}U(\xi,t_*)g(\xi)\,d\xi=\\
&=\lim_{t\to\infty}\mathbf{e}_i^{\top}U(t,t_*)x(t_*)+\int_{t_*}^{\bar h}\lim_{\xi\to\infty}\mathbf{e}_i^{\top}U(\xi,t_*)g(\xi)\,d\xi\\
&=\lim_{t\to\infty}\mathbf{e}_j^{\top}x(t)=\bar x_j,
\end{aligned}
\]
(where all limits exist), proving consensus among $i$ and $j$.
\end{IEEEproof}

\subsubsection{\bf Evolutionary matrices under the AQSC condition}

In the proof of Theorem~\ref{thm.repeated-cont}, we will use additional statements. The first of them relies on Assumption~\ref{asm.int-bound},
and shows that, without loss of generality, one can assume that the sequence $t_p$ in Definition~\ref{def.aqsc} enjoys an arbitrarily large \emph{dwell-time} $\delta>0$,
passing to a subsequence $(t_{pk})$ (where $k\geq 1$ is a fixed number).
\begin{lem}\label{prop.dwell}
Assume that the continuous-time matrix function $A(\cdot)$ is AQSC with some sequence $(t_p)$ and obeys
Assumption~\ref{asm.int-bound}. Then, for any $\delta>0$ there exists an integer $k\geq 1$ such that $t_{p+k}-t_p\geq \delta$.
\end{lem}
\begin{IEEEproof}
Using notation introduced in~\eqref{eq.unions} and~\eqref{eq.int-bound}, all entries of the matrices $A_{t}^{t+\delta}$ do not exceed $\mu_{\delta}$, and their sum does not exceed $n(n-1)\mu_{\bar h}$.
On the other hand, each graph $\g[A_{t_p}^{t_{p+1}}]$ contains a spanning tree that has $n-1$ arcs ($n=\|\V\|$) of weight $\geq\ve$, and hence the total weight of its arcs is not less than $\ve(n-1)$. The union of graphs over $[t_p,t_{p+k}]$, where $k\geq 1$, has the total weight of arcs $\geq k\ve(n-1)$. Choosing $k$ so large that $k\ve>\mu_{\delta} n$, the interval $[t_p,t_{p+k}]$ has therefore length $t_{p+k}-t_p\geq\delta$.
\end{IEEEproof}

Assuming (without loss of generality) that $t_{p+1}-t_p\geq\bar h$ (where $\bar h$ is the delay bound), an
important property of the evolutionary matrices can be derived.
\begin{lem}\label{lem.useful2}
Assume that %Assumption~\ref{asm.int-bound} is valid and
the AQSC condition holds with $t_{p+1}-t_p\geq\bar h$. Then, matrices $U(t_{p+2},t_p)$ have
have uniformly positive diagonal entries $U(t_{p+2},t_p)_{ii}\geq\tilde\eta$ and their graphs are quasi-strongly
$\tilde\ve$-connected, where $\tilde\eta,\tilde\ve>0$ are constants determined by $A(\cdot)$, $\bar h$, $\ve$ from
Definition~\ref{def.aqsc} and $\ell$ from~\eqref{eq.ell}.
%If $\bar h=0$ (undelayed case), then  Assumption~\ref{asm.int-bound} can be discarded.
\end{lem}
\begin{IEEEproof}
For a fixed $j\in\V$, consider the solution $x(t)=U(t_{p+2},t_p)\mathbf{e}_j$ to~\eqref{eq.conse1d}, which corresponds to initial conditions~\eqref{eq.initial} with $t_*=t_p$, $x_*=\mathbf{e}_j$, $\vp(s)\equiv 0\,\forall s<0$. Using the second inequality~\eqref{eq.simple-est} (where $s=t_p$, $\la(s)=0$), one proves that for each $t\in[t_p,t_{p+2}]$
\be\label{eq.aux5}
x_j(t)=U(t,t_p)_{jj}\geq e^{-\int_{t_p}^{t}\alpha_i(\xi)d\xi}\overset{\eqref{eq.ell}}{\geq}\tilde\eta\dfb e^{-2(n-1)\ell}.
\ee
Applying this for $t=t_{p+2}$, one proves the first statement. Also, for any $i\ne j$ and $t\in [t_{p+1},t_{p+2}]$ one has
\[
x_i(t)\geq -\alpha_i(t)x_i(t)+a_{ij}(t)\hat x_j^i(t)\overset{\eqref{eq.aux5}}{\geq}-\alpha_i(t)x_i(t)+a_{ij}(t)\tilde\eta,
\]
because $x_k(t)\geq 0\,\forall k\in\V\,\forall t\geq t_p$ and $\hat x^i_j(t)=x_j(t-h_{ij}(t))$, where $t-h_{ij}(t)\geq t_{p+1}-\bar h\geq t_p$. Taking into account that $x_i(t_{p+1})\geq 0$, one proves that
\[
\begin{gathered}
x_i(t_{p+2})\geq\int_{t_{p+1}}^{t_{p+2}}e^{-\int_{t_{p+1}}^s\alpha_i(\xi)d\xi}a_{ij}(s)\hat x_j^i(s)ds\geq\\
\geq e^{-(n-1)\ell}\tilde\eta\int_{t_{p+1}}^{t_{p+2}}a_{ij}(s)ds.
\end{gathered}
\]
In particular, if $(A_{t_{p+1}}^{t_{p+2}})_{ij}\geq\ve$, then one has
\[
x_i(t_{p+2})=U(t_{p+2},t_p)_{ij}\geq\tilde\ve\dfb e^{-(n-1)\ell}\tilde\eta=e^{-3(n-1)\ell}.
\]
By definition of the AQSC matrix $A(\cdot)$, the graph $\g[U(t_{p+2},t_p)]$ is quasi-strongly $\tilde\ve$-connected.
\end{IEEEproof}

\subsubsection{\bf Evolutionary matrices under the NITS condition}

The proof of Theorem~\ref{thm.symm-c}  is based on another property of the evolutionary matrices, specific for type-symmetric graphs.
\begin{lem}\label{lem.tech}
If the continuous-time matrix $A(\cdot)$ obeys the NITS property and Assumption~\ref{asm.int-bound}, then matrices $U(t,t_*)$ have uniformly positive
diagonal elements, that is,
\be\label{eq.diag-pos}
U(t,t_*)_{ii}\geq \varrho>0\quad\forall t\ge t_*\geq 0\,\forall i\in\V.
\ee
The constant $\varrho=\varrho(n,\ell,\mu_{\bar h},K)$ does not depend on a specific matrix $A(\cdot)$ and delays $h_{ij}$,
being determined by the number of agents $n=|\V|$, constant $K$ from~\eqref{eq.type-symm-non}, the number $\ell$
from~\eqref{eq.ell} and number $\mu_{\bar h}$ from~\eqref{eq.int-bound} (where $D=\bar h$).
\end{lem}
\begin{IEEEproof}
Throughout the proof, the sequence $(t_p)$ is same as in Definition~\ref{def.symm}.
Without loss of generality, assume that $t_0=0$.

\textbf{Step 1.} Notice that~\eqref{eq.diag-pos} is entailed by the following:

\emph{(A) A constant $\tilde\varrho=\tilde\varrho(n,\ell,\mu_{\bar h},K)\in (0,1)$ exists featured by the following property. Let $x(t),t\geq t_p$ be such a solution of~\eqref{eq.conse1d} that $\la(t_p)\geq 0$ and $x_i(t)\geq 1$ for $t\in[t_p-\bar h,t_p]$ for some agent $i\in\V$. Then,
function $x_i(t)$ satisfies the inequality
\be\label{eq.column+}
x_i(t)\geq\tilde\varrho\quad\forall t\geq t_p.
\ee
}

Indeed, suppose that statement (A) is valid. Consider solution $x(t)=U(t,t_*)\mathbf{e}_i$, where $i\in\V$ and $t_*\geq 0$. Let $p\geq 1$ be the minimal index such that $t_p\geq t_*+\bar h$. In view of~\eqref{eq.int-bound} (with $D=\bar h$)
and~\eqref{eq.ell}, one has
\[
\int_{t_*}^{t_p}a_{ij}(t)\,dt\leq\ell+\mu_{\bar h}\quad\forall i\ne j.
\]
Using inequalities~\eqref{eq.simple-est} and~\eqref{eq.ell}, one has
\[
x_i(t)\overset{\eqref{eq.simple-est}}{\geq} e^{-\int_{t_*}^{t_p}\alpha_i(s)\,ds}x_i(t_*)\geq \theta\dfb e^{-(n-1)(\ell+\mu_{\bar h})}\,\forall t\in[t_*,t_p],
\]
in particular, $x_i(t)\geq \theta$ for $t\in[t_p-\bar h,t_p]$. Applying (A) to solution $\tilde x(t)=\theta^{-1}x(t)$, one shows that
\[
x_i(t)=U(t,t_*)_{ii}\geq\varrho\dfb \theta\tilde\varrho\quad\forall t\ge t_*.
\]

\textbf{Step 2.} We will show that statement (A), in turn, is implied by another condition (B) presented below.
Given a solution to system~\eqref{eq.conse1d}, $\V_0\subseteq\V$ and $t\geq 0$, denote
\[
\la_{\V_0}(t)\dfb\min_{i\in\V_0}\inf_{s\in[t-\bar h,t]} x_i(s).
\]

\emph{(B) A constant $\gamma=\gamma(n,\ell,\mu_{\bar h},K)\in (0,1/2)$ exists featured by the following property.
If $x(t),\,t\geq t_p$ is a solution of~\eqref{eq.conse1d}
such that $\la(t_p)\geq 0$ and $\la_{\V_0}(t_p)\geq 1\,\forall i\in\V_0\subseteq\V$ for some $p\geq 0$, then one of the
following statements (i) and (ii) hold
\begin{enumerate}[\textbf(i)]
  \item $\la_{\V_0}(t)\geq 1/2\,\forall t\geq t_p$;
  \item a set $\V_1\supsetneq\V_0$ and index $q>p$ exist such that
  $\la_{\V_0}(t)\geq \gamma\,\forall j\in\V_0$ when $t\in[t_p,t_q]$ and $\la_{\V_1}(t_q)\geq \gamma$.
\end{enumerate}}

Indeed, assume that (B) holds and consider a solution satisfying the assumptions of statement (A). Applying (B) to this
solution and set $\V_0=\{j:x_j(t_p)\geq 1\}\ni\{i\}$, one of statements \emph{(i)} and \emph{(ii)} should hold.

If \emph{(i)} holds, then $x_i(t)\geq \beta_0\dfb 1/2$ for all $t\geq t_p$.

Assume that \emph{(ii)} holds and let $q_1>p$, $\V_1\supsetneq\V_0$ be the corresponding integer index and set of agents.
Applying (B) to the solution $\tilde x(t)=\gamma^{-1}x(t)$, $\tilde p=q_1$ and $\tilde{\V}_0=\V_1$, one shows that
either condition \emph{(i)} holds, and then $\la_{\V_1}(t)\geq\beta_1\dfb \beta_0\gamma$ for $t\geq t_{q_1}$ or the scenario
from \emph{(ii)} is realized, and there exist such a set $\V_2\supsetneq\V_1$ and an index $q_2>q_1$
that $\la_{\V_1}(t)\geq\gamma^2$
when $t\in[t_{q_1},t_{q_2}]$ and $\la_{\V_2}(t_{q_2})\geq \gamma^2$.

In the latter situation, we repeat the procedure and apply (B) to the solution $\hat x(t)=\gamma^{-2}x(t)$,
$\hat p=q_2$ and $\hat V_0=V_2$, showing that either $\la_{\V_2}(t)\geq\beta_2\dfb \beta_0\gamma^2$ for
$t\geq t_{q_2}$ or there exist a set $\V_3\supsetneq\V_2$ and index $q_3>q_2$ such that $\la_{\V_2}(t)\geq\gamma^3$
when $t\in[t_{q_2},t_{q_3}]$ and $\la_{\V_3}(t_{q_3})\geq \gamma^3$, in which situation we can again apply statement (B),
and so on.

Since $|\V|=n$, this procedure terminates after $1\leq r\leq n-1$ steps, after which
scenario \emph{(i)} is realized and the set $\V_r$ cannot be constructed. By construction, 
$i\in\V_0\subset\ldots\subset\V_{r-1}$, and hence
$x_i(t)\geq\beta_{r-1}=\beta_0(\gamma)^{r-1}\geq \beta_0(\gamma)^{n-2}=\gamma^{n-2}/2$ for all $t\geq t_p$.
Hence, (B) implies (A) with $\tilde{\varrho}=\gamma^{n-2}/2$.

\textbf{Step 3.} We are now going to prove statement (B) via induction on $|\V|$.
The induction base $|\V|=1$ is trivial, in this situation the only agent obeys the equation $\dot x=0$, 
so if $x(t_p)\geq 1$, then condition \emph{(i)} holds automatically.

Suppose that (B) (and thus also statement of Lemma~\ref{lem.tech}) has been proved for groups of $\leq n-1$ agents and $|\V|=n$.
Notice that for any subgroup $\tilde\V\subsetneq\V$, the corresponding matrix $\tilde A=(a_{ij})_{i,j\in\tilde V}$ obeys
the NITS conditions and~\eqref{eq.int-bound} with the same constant $K$, sequence $(t_p)$ and constants $\ell,\mu_{\bar h}$
as $A$. We know that Lemma~\ref{lem.tech} is valid for each reduced system
\begin{gather}
\dot x_i=\sum\limits_{j\in\tilde\V,j\ne i}a_{ij}(t)(\hat x_j^i(t)-x_i(t)),\,i\in\tilde\V.\label{eq.subsys0}
\end{gather}
Introducing the corresponding evolutionary matrix $\tilde U(t,t_*)$, we
have $\tilde U(t,t_*)_{ii}\geq\varrho(|\V'|,\ell,\mu_{\bar h},K)$. We define
\[
\varrho'=\rho'(\ell,\mu_{\bar h},K)\dfb\min_{k\leq n-1}\varrho(k,\ell,\mu_{\bar h},K)>0.
\]

Consider first the solution with the special initial condition:
\be\label{eq.init-special}
\begin{gathered}
x_i(t)\equiv 1\,\forall t\in[t_p-\bar h,t_p]\,\forall i\in\V_0\\
x_i(t)\equiv 0\,\forall t\in[t_p-\bar h,t_p]\,\forall i\not\in\V_0.
\end{gathered}
\ee

Denoting $\V^+\dfb\V_0$ and $\V^{\dagger}=\V\setminus\V^+$, consider systems
\begin{gather}
\dot x_i=\sum\limits_{j\in\V^+}a_{ij}(t)(\hat x_j^i(t)-x_i(t))+f_i^+(t),\,i\in\V^+;\label{eq.subsys1}\\
\dot x_i=\sum\limits_{j\in\V^{\dagger}}a_{ij}(t)(\hat x_j^i(t)-x_i(t))+f_i^{\dagger}(t),\,i\in\V^{\dagger}.\label{eq.subsys2}
\end{gather}
where $f^+(t)\in\r^{\V^+}$ and $f^{\dagger}(t)\in\r^{\V^{\dagger}}$ are some functions. 
Let $U^+,U^{\dagger}$ stand for the evolutionary matrices of systems~\eqref{eq.subsys1} and~\eqref{eq.subsys2} respectively. 

For the solution $x(t)$ determined by~\eqref{eq.init-special}, two situations are possible:
1) $\la_{\V^+}(t)=\la_{\V_0}(t)\geq 1/2\;\;\forall t\geq t_p$ or 2) the minimal $t'>t_p$ exists
such that $\min_{i\in\V^+}x_i(t')=1/2$.

In the second case, let $s\geq p+1$ be such an index that $t_{s-1}<t'\leq t_s$. Denote
\be\label{eq.const-aux}
\begin{gathered}
c_1\dfb e^{-(n-1)\ell}/2,\;c_2\dfb\frac{\varrho'Cc_1}{\varrho'C+1},\,
C\dfb\frac{2}{Kn^2}.
\end{gathered}
\ee
where $K\geq 1$ is the constant from~\eqref{eq.type-symm-non}. By assumption, $x_i(t')\geq 1/2\,\forall i\in\V^+$, and thus, using the second inequality in~\eqref{eq.simple-est},~\eqref{eq.ell} and~\eqref{eq.init-special}, one has
$\la_{\V^+}(t)\geq c_1\forall t\in[t_p,t_q]$.

We are now going to show that an index $j\in\V^{\dagger}$ and an instant $t\in [t_p,t_s]$ exist such that $x_j(t)\geq c_2$.
Assume, on the contrary, that $x_j(t)<c_2$ on $[t_p,t_s]$ for all $j\in\V^{\dagger}$.
Subvectors $x^+(t)=(x_i(t))_{i\in\V^+}$, $x^{\dagger}(t)=(x_i(t))_{i\in\V^{\dagger}}$ are solutions, respectively, to~\eqref{eq.subsys1} and~\eqref{eq.subsys2}, where
\be\label{eq.aux-f}
\begin{gathered}
f_i(t)=
\begin{cases}
\sum_{j\in\V^{\dagger}}a_{ij}(t)[\hat x^i_j(t)-x_i(t)],\; i\in\V^+\\
\sum_{j\in\V^{+}}a_{ij}(t)[\hat x^i_j(t)-x_i(t)],\; i\in\V^{\dagger}.
\end{cases}
\end{gathered}
\ee
In view of~\eqref{eq.init-special} and Lemma~\ref{prop.bound}, $0\leq x_k(t)<c_2<c_1\leq x_m(t)\leq 1$ for any
$k\in\V^{\dagger}$, $m\in\V^+$ and $t\in [t_p,t_s]$. Hence,
\[
\begin{gathered}
f_i(t)\geq -\sum_{j\in\V^{\dagger}}a_{ij}(t)\quad\forall i\in \V^+,\\
f_j(t)\geq (c_1-c_2)\sum_{i\in\V^{+}}a_{ji}(t)\geq 0\quad\forall j\in \V^{\dagger},
\end{gathered}
\]
whenever $t\in[t_p,t_q]$.  Using the Cauchy formula~\eqref{eq.cauchy-del-estim} (where $t_*$ is replaced by $t_p$, functions $x,f$ are replaced by respectively $x^+,f^+$ and $\la(t_p)=1$), one has
\[
\begin{aligned}
x^+(t')\geq \ones_{\V^+}+\int_{t_p}^{t'}U^+(t',t)f^+(t)dt\Longrightarrow\\
x_i(t')\geq 1-\sum_{k\in\V^+,j\in\V^{\dagger}}\int_{t_p}^{t_q}a_{kj}(t)dt\quad\forall i\in\V^+.
\end{aligned}
\]
By assumption, $x_i(t')=1/2$ for some $i$. Taking into account that $r\dfb\|V^{\dagger}\|\leq n-1$ and $r(n-r)\leq n^2/4$,
one shows that for some $k\in\V^+$ and $j\in\V^{\dagger}$ we have
\be\label{eq.auxaux}
\int\limits_{t_p}^{t_q}a_{kj}(t)dt\geq \frac{1/2}{n^2/4}=\frac{2}{n^2}\overset{\eqref{eq.type-symm-non}}{\Longrightarrow}
\int\limits_{t_p}^{t_q}a_{jk}(t)dt\geq \frac{2}{Kn^2}=C.
\ee

Using the Cauchy formula~\eqref{eq.cauchy-del-estim} (where $t_*$ is replaced by $t_p$, functions $x,f$ are replaced by respectively $x^{\dagger},f^{\dagger}$ and $\la(t_p)=0$), one, on the other hand, obtains the inequality
\[
\begin{aligned}
x^{\dagger}(t_q)\geq \int_{t_p}^{t_q}U^{\dagger}(t_q,t)f^{\dagger}(t)dt
\end{aligned}
\]
Recalling that $U^{\dagger}(t,s)_{jj}\geq\rho'$ because $\V^{\dagger}\subsetneq\V$, one has
\[
\begin{aligned}
x_j(t_q)\geq \varrho'\int_{t_p}^{t_q}f_j(t)dt\geq\varrho'(c_1-c_2)\int_{t_p}^{t_q}a_{jk}(t)dt\geq\\
\geq\varrho'(c_1-c_2)C=c_2,
\end{aligned}
\]
which leads to the contradiction.

Hence, there exists such an instant $t''\in[t_p,t_s]$ and $j\in\V^{\dagger}$ such that $x_j(t'')\geq c_2$.
Recall also that $\la_{\V^+}(t'')\geq c_1>c_2$. Denoting $\V_1=\V_0\cup\{j\}\supsetneq\V_0$ and choosing the minimal index $q\geq p$ such that $t_q\geq t+\bar h$, one has
one has
\[
\int_{t}^{t_q}a_{km}(s)\,ds\leq\ell+\mu_{\bar h}\quad\forall k\ne m,
\]
and hence, using~\eqref{eq.simple-est} and~\eqref{eq.ell}, one can obtain that
\be\label{eq.auxaux1}
\begin{gathered}
\la_{\V_1}(t_q)\geq \gamma\dfb c_2e^{-(n-1)(\ell+\mu_{\bar h})},\\
\la_{\V_0}(t)\geq \gamma\quad\forall t\in[t_p,t_q].
\end{gathered}
\ee

We have shown that statement (B) holds for solutions with special initial conditions~\eqref{eq.init-special};
note for such a solution $\V_1$ and $t_q$ depend on the choice of $\V_0\subseteq\V$, but the constant $\gamma>0$ does 
not depend on it. Consider now any other solution $\tilde x(t)$ such that $\la_{\V_0}(t_p)\geq 1$ and $\la(t_p)\geq 1$.
Then, in view of Remark~\ref{rem.monotone}, one has $\tilde x(t)\geq x(t)\,\forall t\geq t_p$, where $x(t)$ is determined by the 
initial condition~\eqref{eq.init-special}, and hence $\tilde x(t)$ also satisfies one of the conditions \emph{(i),(ii)}.
This finishes the proof of induction step. We have finished the proof of statement (B), implying also statement (A)
and Lemma~\ref{lem.tech}.
\end{IEEEproof}

\subsection{Proof of Lemma~\ref{lem-l1-rob}}\label{subsec.proof-robust}

Assume that algorithm~\eqref{eq.conse1d} establishes consensus between agents $i$ and $j$.%., that is,~\eqref{eq.conse-cauchy} holds.
In view of~\eqref{eq.cauchy+}, every solution to the disturbed system~\eqref{eq.conse1d-f} (started at time $t_*\geq 0$) is the sum of a solution to~\eqref{eq.conse1d} and the ``forced'' solution
\be\label{eq.forced}
x^f(t)=\int_{t_*}^tU(t,\xi)f(\xi)\,d\xi=\int_{t_*}^{\infty}\hat U(t,\xi)f(\xi)\,d\xi,
\ee
where $\hat U(t,\xi)=0$ if $\xi>t$ and $\hat U(t,\xi)=U(t,\xi)$ for $\xi\leq t$.

To prove that consensus between $i$ and $j$ is preserved when $f\in L_1([0,\infty)\to\r^{\V})$, it thus suffices to show that
\[
\lim_{t\to\infty} x^f_i(t)=\lim_{t\to\infty} x^f_j(t)
\]
(and both limits exist). To prove this,  notice that~\eqref{eq.conse-cauchy} can be rewritten as follows: for all $t_*\geq 0$, the vector $u^{\top}_{t_*}$ exists such that
$\lim_{t\to\infty}\mathbf{e}_i^{\top}\hat U(t,t_*)=\lim_{t\to\infty}\mathbf{e}_j^{\top}\hat U(t,t_*)=u^{\top}_{t_*}$. Recalling that matrices $\hat U(t,\xi)$ are uniformly bounded (Lemma~\ref{lem.cauchy0}) and $f$ is $L_1$-summable, the Lebesgue dominated convergence theorem entails
the existence of coincident limits
\[
\begin{gathered}
\lim_{t\to\infty}x_i^f(t)=\lim_{t\to\infty}\int_{t_*}^{\infty}\mathbf{e}_i^{\top}\hat U(t,t_*)f(\xi)\,d\xi=\int_{t_*}^{\infty}u^{\top}_{t_*}f(\xi)d\xi,\\
\lim_{t\to\infty}x_j^f(t)=\lim_{t\to\infty}\int_{t_*}^{\infty}\mathbf{e}_j^{\top}\hat U(t,t_*)f(\xi)\,d\xi=\int_{t_*}^{\infty}u^{\top}_{t_*}f(\xi)d\xi,
\end{gathered}
\]
which finishes the proof of Lemma~\ref{lem-l1-rob}. $\blacksquare$

\subsection{Proof of Theorem~\ref{thm.repeated-cont} and Lemma~\ref{lem-rob-general}}\label{subsec.proof-repeated}

Theorem~\ref{thm.repeated-cont} is based on Lemma~\ref{lem.cauchy0} and the following lemma.

\begin{lem}\label{lem.tau1-ext}
Let $B^1,\ldots,B^{n-1}\in\r^{\V\times\V}$, where $n=|\V|-1$, be \emph{substochastic} matrices with positive diagonal
entries and quasi-strongly connected graphs $\g[B^i]$. Consider a sequence $z^1,\ldots,z^n\in\r^v$, where
$0\leq z^1\leq\ones$ and
\be\label{eq.aux33}
0\leq z^k\leq B^{k-1}z^{k-1}+(\ones-B^{k-1}\ones),\quad k=2,\ldots,n.
\ee
Then, $0\leq\max z^n-\min z^n<1$. Moreover, for each $\ve,\eta>0$
there exists $\rho=\rho(\ve,\eta)\in (0,1)$ such that if $\g[B^i]$ are quasi-strongly $\ve$-connected and
$b^{1}_{ii},\ldots,b^{n-1}_{ii}\geq\eta\,\forall i\in\V$, then
$\max z^n-\min z^n\leq\rho(\eta,\ve)$.
\end{lem}
\begin{IEEEproof} For a vector $x\in[0,1]^{\V}\dfb\{(x_i)_{i\in\V}:x_i\in[0,1]\,\forall i\}$, denote $\textbf{ZER}(x)\dfb\{i:x_i=0\}\subseteq\V$ and
$\textbf{ONE}(x)\dfb\{i:x_i=1\}\subseteq\V$. In view of~\eqref{eq.aux33} and substochasticity of $B^k$,
we have $z^k\in [0,1]^{\V}$ for all $k=1,\ldots,n$.

For any substochastic matrix $B$ with positive diagonal entries and any $x\in [0,1]^{\V}$,
the inequality $x_i>0$ (respectively, $x_i<1$) entails that $(Bx)_i>0$ (respectively, $(Bx)_i<1$). Hence, the sets
$\textbf{ZER}(z^k)$ and $\textbf{ONE}(z^k)$ are nested: $\textbf{ZER}(z^k)\subseteq \textbf{ZER}(z^{k-1})$ and
$\textbf{ONE}(z^k)\subseteq \textbf{ONE}(z^{k-1})$.

Furthermore, as can be easily seen, if $j\not\in \textbf{ONE}(z^k)$ and $b_{ij}^{k}>0$,
then $z_i^{k+1}<1$. Similarly, if $j\not\in \textbf{ZER}(z^k)$ and $b_{ij}^k>0$, then $z_i^{k+1}>0$. If $\g[A^k]$ is quasi-strongly connected,
with some node $r$, then either $r\not\in\textbf{ONE}(z^k)$ or $r\not\in\textbf{ZER}(z^k)$. In the first situation,
we either have $\textbf{ONE}(z^k)=\emptyset$ or a path connecting $r$ to $\textbf{ONE}(z^k)$ in $\g[B^k]$ should exist,
that is, at least one arc comes to $\textbf{ONE}(z^k)$ from outside. In this situation,
$\textbf{ONE}(z^{k+1})\subsetneq \textbf{ONE}(z^k)$. Similarly, if $r\not\in\textbf{ZER}(z^k)$, then either
$\textbf{ZER}(z^k)=\emptyset$ or
$\textbf{ZER}(z^{k+1})\subsetneq \textbf{ZER}(z^k)$. Since the cardinality of $\textbf{ONE}(z^1)\cup \textbf{ZER}(z^1)$
is not greater than $n=|\V|$, at least one of the sets $\textbf{ONE}(z^n)$ or $\textbf{ZER}(z^n)$ is empty,
and hence $\max z^n-\min z^n>0$.

The set $\mathfrak{B}(\eta,\ve)$ of substochastic matrices $B$ such that $a_{ii}\geq\eta$ and $\g[B]$ is quasi-strongly
$\ve$-connected is \emph{compact}. Hence, the set of all sequences $(z^1,\ldots,z^n)$ obeying~\eqref{eq.aux33} is also
compact. The continuous function $\max z^n-\min z^n$ thus reaches a maximum $\rho=\rho(\eta,\ve)<1$ on $\mathfrak{B}(\eta,\ve)$.
\end{IEEEproof}
\vskip0.1mm

\subsubsection{Proof of Theorem~\ref{thm.repeated-cont}}
In view of Lemma~\ref{prop.dwell}, we may also assume that
$t_{p+1}-t_p\geq\bar h$ and therefore Lemma~\ref{lem.useful2} is applicable.
Without loss of generality, assume that $r=0$.
For brevity, we denote $D(t)\dfb\La(t)-\la(t)$. In view of Lemma~\ref{prop.bound}, $D$ is a non-increasing function,
in particular, $D(t_0)\leq D(t_*)$.

Notice also that, without loss of generality, we may assume that $\la(t_0)=0$ and $\La (t_0)=1$. Indeed, if
$\la(t_0)=\La(t_0)$, then consensus obviously established and~\eqref{eq.contract} holds, because
$La(t)\equiv\La(t_0)=\la(t)=x_i(t)\,\forall i$ due to Lemma~\ref{prop.bound}. Otherwise,
the transformation $x(t)\mapsto x^0(t)\dfb[x_i(t)-\la(t_0)\ones]/[\La(t_0)-\la(t_0)]$ gives a solution to~\eqref{eq.conse1d},
which corresponds to $\La^0(t_0)=1$ and $\la^0(t_0)=0$.

If $\la(t_0)=0$ and $\La (t_0)=1$, then $0\leq\la(t)\leq\bar La(t)\leq 1\,\forall t\geq t_0$ thanks to Lemma~\ref{prop.bound}.
Using~\eqref{eq.cauchy-del-estim} (where $t_*$ is replaced by $t_p$ and $t=t_{p+2}$, one arrives at
\[
0\leq x(t_{p+2})-U(t_{p+2},t_p)x(t_p)\leq \ones-U(t_{p+2},t_p)\ones.
\]
In view of Lemma~\ref{lem.useful2}, matrices $B^j=U(t_{2j},t_{2j-2})$ have strongly positive diagonal entries and
quasi-strongly $\tilde\ve$-connected graphs (where $\tilde\ve>0$ is some constant determined by $A(\cdot)$ and $\bar h$).
Applying Lemma~\ref{lem.tau1-ext} to matrices $B^j=U(t_{2j},t_{2j-2})$ and vectors $z^j=x(t_{2j-2})$, one arrives at
\[
\max x(t_{2(n-1)})-\min x(t_{2(n-1)})\leq\rho<1,
\]
where $\rho$ is a constant (independent of the specific solution). Substituting $t_*=t_{2n-2}$
into~\eqref{eq.cauchy-del-estim}, we have
\[
0\leq x(t)-U(t,t_{2n-2})x(t_{2n-2})\leq \ones-U(t,t_{2n-2})\ones,
\]
and thus for $t\in [t_{2n-2},t_{2n-1}]$ inequalities hold
\[
\begin{aligned}
U(t,t_{2n-2})_{ii}&x_i(t_{2n-2})\leq x_i(t)
\leq
\\
&\leq U(t,t_{2n-2})_{ii}\,x_i(t_{2n-2})+1-U(t,t_{2n-2})_{ii}.
\end{aligned}
\]
Since the diagonal entries $U(t,t_{2n-2})_{ii}$ are uniformly positive over $t\in [t_{2n-2},t_{2n-1}]$
(Lemma~\ref{lem.useful2}), $\max x(t)-\min x(t)\leq\theta<1$ on
the interval  $t\in [t_{2n-2},t_{2n-1}]$, whence
\[
D(t_{2n-1})=\La(t_{2n-1})-\la(t_{2n-1})\leq\theta=\theta D(t_0).
\]
Here $\theta$ is some constant, which depends only on the matrices $U(t,s)$ and does not depend on $t_0$ or the specific
solution. Recalling that $0\leq\la(t_{2n-1})\leq\La(t_{2n-1})\leq 1$,
we can now repeat the argument, replacing $t_0$ by $t_{2n-1}$ and showing that
\be\label{eq.aux34}
\begin{split}
D(t_{(2n-1)k})\leq\theta D(t_{(2n-1)(k-1)})\leq\ldots\leq\theta^kD(t_0),
\end{split}
\ee
which finishes the proof of the second inequality in~\eqref{eq.contract}. The first inequality is obvious from~\eqref{eq.barx-i}
and inequalities $\la(t)\ones\leq x(t)\leq\La(t)\ones$.
Finally, note that Assumption~\ref{asm.int-bound} was used
only to invoke Lemma~\ref{prop.dwell}. If $\bar h=0$, then $t_{p+1}-t_p>\bar h$ automatically and hence Assumption~\ref{asm.int-bound}
can be discarded. $\blacksquare$

\subsubsection{Proof of Lemma~\ref{lem-rob-general}}

Applying~\eqref{eq.contract} to $t_*=s$ and solution $x(t)=U(t,s)v$, where $v\in\r^{\V}$ is an arbitrary vector with $\|v\|_{\infty}=1$
(for such a solution, obviously, $\La(s)-\la(s)\leq 2$), one proves that
\be\label{eq.u-converge}
\begin{gathered}
\|U(t,s)-\bar U_{s}\|_{\infty}\leq 2\theta^k,\\
t\in[t_{r+2k(n-1)},t_{r+2(k+1)(n-1)}],\,t_r\geq s.
\end{gathered}
\ee
Note that $\theta$ does not depend on $s$, and $s$ can be arbitrary.

For an arbitrary vector $x\in\r^{\V}$, denote $\delta(x)\dfb\max x-\min x$; obviously, $0\leq\Delta(x)\leq 2\|x\|_{\infty}$ and
$\delta(x+y)\leq\delta(x)+\delta(y)$. Similar to Lemma~\ref{lem-l1-rob}, it suffices to prove~\eqref{eq.rob-general} for
the ``forced'' solution~\eqref{eq.forced}. To this end, we introduce the function
\[
x_*^f(t)\dfb\int_{t_*}^t{\bar U}_{s}f(s)\,ds.
\]
It seems natural that, in view of~\eqref{eq.u-bar}, $x^f$ and $x_*^f$ are sufficiently close when $t$ becomes large.
To make the latter statement formal, choose an index $q\geq r$ (so that $t_q\geq t_*$). Obviously,
\[
\begin{aligned}
\left\|x^f(t)-x_*^f(t)\right\|_{\infty}\leq\int_{t_*}^{t_q}\|U(t,s)-\bar U_{s}\|_{\infty}\,\|f(s)\|_{\infty}ds+\\
+\int_{t_q}^{t}\|U(t,s)-{\bar U}_{s}\|_{\infty}\,\|f(s)\|_{\infty}dt
\end{aligned}
\]
The Lebesgue dominated convergence theorem ensures that the first integral vanishes as $t\to\infty$.
In view of~\eqref{eq.u-converge}, the second
integral can be estimated as
\[
\begin{aligned}
\sum_{k=0}^{\infty}
\underbrace{2(n-1)\,2(1+\theta+\theta^2+\ldots)}_{=C_1}\sup_{p\geq q}\int_{t_p}^{t_{p+1}}\|f(t)\|_{\infty}dt.
\end{aligned}
\]
One thus concludes that for every $q\geq r$ we have
\[
\begin{aligned}
\varlimsup_{t\to\infty}\left\|x^f(t)-x_*^f(t)\right\|_{\infty}\leq C_1\sup_{p\geq q}\int_{t_p}^{t_{p+1}}\|f(t)\|_{\infty}dt.
\end{aligned}
\]
Since $q$ can be arbitrary, one can now pass to the limit $q\to\infty$, replacing $\sup$ in the latter inequality by
$\varlimsup_{p\to\infty}$. Hence,
\[
\begin{aligned}
&\varlimsup_{t\to\infty}\delta(x^f(t))\leq\varlimsup_{t\to\infty}\underbrace{\delta(x^f_*(t))}_{=0}
+\varlimsup_{t\to\infty}\delta(x^f(t)-x^f_*(t))\leq\\
&\leq 2\varlimsup_{t\to\infty}\left\|x^f(t)-x_*^f(t)\right\|_{\infty}=
2C_1\varlimsup_{p\to\infty}\int_{t_p}^{t_{p+1}}\|f(t)\|_{\infty}dt,
\end{aligned}
\]
which proves~\eqref{eq.rob-general} for the forced solution $x^f(t)$ with $C=2C_1$. Since the undisturbed solutions reach global consensus,
\eqref{eq.rob-general} also holds for all solutions to~\eqref{eq.conse1d-f}.

\subsection{Proof of Theorem~\ref{thm.symm-c}}

In this subsection, we assume that the conditions of Theorem~\ref{thm.symm-c} hold.
We will prove, in fact, that consensus is guaranteed for every bounded solution to the inequality~\eqref{eq.conse1d-ineq} (in particular, for all solutions to~\eqref{eq.conse1d}).
For every such a solution, the function $\La(t)$ is non increasing (Lemma~\ref{prop.bound}), and thus the limit $\bar\La=\lim_{t\to\infty}\La(t)>-\infty$ exists.
Assume that $\La(t_*)>\la(t_*)$ (otherwise, consensus is obvious due to Lemma~\ref{prop.bound}).

Proof of Theorem~\ref{thm.symm-c} is based on the fruitful idea of a one-to-one \emph{ordering mapping} $\sigma_t:[1:n]\to\V$ ($n=|\V|$), which sorts the vector is the ascending order
\[
z_1(t)=x_{\sigma_t(1)}(t)\leq\ldots\leq z_n(t)=x_{\sigma_t(n)}(t).
\]

Notice that, in general, $\sigma_t(i)$ is defined non-uniquely (if $x(t)$ has two or more equal components), however, the mapping $\sigma_t$ can always be chosen \emph{measurable} in $t$ and the functions $z_i(t)$ are absolutely continuous on $[0,\infty)$~\cite{TsiTsi:13,ProMatvCao:2016}.

\subsubsection{The case of \textbf{connected} persistent graph $\g_{\infty}$}

We first consider the case where persistent graph $\g_{\infty}$ is connected and prove that \emph{global} consensus is established in this case.
Our goal is to show that $x(t)\xrightarrow[t\to\infty]{}\bar\La\ones$. %Since matrix $A(\cdot)$ is $(t_p)$-bounded~\eqref{eq.ell},
It suffices to prove that
\be\label{eq.x-consensus-tp}
x(t_p)\xrightarrow[p\to\infty]{}\bar\La\ones.
\ee
Indeed, if~\eqref{eq.x-consensus-tp} is valid, then, choosing $t_{p-1}<t'<t_p$ and using Lemma~\ref{lem.cauchy0} (where $t_*$ in ~\eqref{eq.cauchy-del-estim} is replaced by $t'$), one has
%\ben%\label{eq.aux2}
%\begin{gathered}
%\int_{t}^{t_p}\alpha_i(s)ds\overset{\eqref{eq.ell}}{\leq}(n-1)\ell\overset{\eqref{eq.simple-est}}{\Longrightarrow}
%x_i(t_p)\leq (1-\theta_0)\La(t)+\theta_0x(t)\\ \theta\dfb e^{-(n-1)\ell}\,\forall i\in\V_*,
$x(t_p)\leq \psi x(t')+(1-\psi)\La(t')\ones$,
%\end{gathered}
%\een
and hence $[x(t_p)-(1-\psi)\La(t')\ones]\psi^{-1}\leq x(t')\leq\La(t')\ones$, where the leftmost and rightmost expressions both go to $\bar\La\ones$ as $p\to\infty$.

Recalling the definition of $z_k$,~\eqref{eq.x-consensus-tp} is equivalent to
\be\label{eq.z-consensus}
z_k(t_p)\xrightarrow[p\to\infty]{}\La^*\quad\forall k=1,\ldots,n.
\ee
We are going to prove~\eqref{eq.z-consensus} using backward induction on $k=n,\ldots,1$. The induction base ($k=n$) follows from Lemma~\ref{prop.bound}.

Suppose that~\eqref{eq.z-consensus} has been proved for $k=r+1,\ldots,n$. To prove it for $k=r$, it suffices to show that
\be\label{eq.z-consensus1}
\varliminf_{p\to\infty}z_r(t_p)\geq\La_*.
\ee
Assume, on the contrary, that a subsequence $\tau_m\dfb t_{p(m)}$ exists, $m\to\infty$ such that $z_r(\tau_m)\leq\La(\tau_m)-\delta\,\forall m$.
Passing to a subsequence, one can assume, without loss of generality, that $\V^+\dfb\sigma_{\tau_m}([1:r])\subset\V$ does not depend on $m$
(this set contains indices of the $r$ smallest values at time $\tau_m$).  By definition, $x_i(\tau_m)\leq z_r(\tau_m)\leq\La(\tau_m)-\delta\;\;\forall i\in\V^+\,\forall m$.
Without loss of generality, we choose $m$ so large that
\be\label{eq.aux-assum}
z_{r+1}(t)>\La(\tau_m)-\delta/3\quad\forall t\geq\tau_m-\bar h.
\ee

We will show that the aforementioned assumptions contradict to the induction hypothesis, proving  the existence of such $j\not\in\V^+$,
a sequence $\tilde\tau_m=t_{\tilde p(m)}\to\infty$ and $\tilde\delta<\delta$ that
\be\label{eq.aux++}
\max_{i\in\V^+\cup\{j\}}x_i(\tilde\tau_m)\leq\La(\tau_m)-\tilde\delta.
\ee
The cardinality of $\V^+\cup\{j\}$ is $r+1$, so~\eqref{eq.aux++} entails that $z_{r+1}(t_{\tilde p(m)})\leq\La(\tau_m)-\tilde\delta$ for each $m$, which is incompatible with the induction hypothesis.
Informally, speaking, as a result of persistent interaction between $\V^+$ and $\V^{\dagger}$, some agent from $\V^{\dagger}$ comes ``sufficiently close'' to agents from $\V^+$.

The proof of~\eqref{eq.aux++} is based on Lemma~\ref{lem.tech}, and, in fact, similar to Step 3 of its proof. We divide it into several steps.
\textbf{Step 1.} Denote $\V^{\dagger}\dfb\V\setminus\V^+$.
Obviously, the submatrix $A^{\dagger}=(a_{ij})_{i,j\in\V^{\dagger}}$ (as well as the submatrix $A^{+}=(a_{ij})_{i,j\in\V^{+}}$) obeys the NITS condition, as well as Assumption~\ref{asm.int-bound}.
Let $U^{\dagger}(t,s)$ be the evolutionary matrix, corresponding to submatrix $A^{\dagger}$ or, equivalently, the ``reduced'' delay system
\begin{gather}
\dot x_i=\sum\limits_{j\in\V^{\dagger}}a_{ij}(t)(\hat x_j^i(t)-x_i(t)),\,i\in\V^{\dagger}\label{eq.subsys2}.
\end{gather}
Matrix $U^+(t,s)$ is defined similarly and corresponds to submatrix $A^+(t)$. Both submatrices $A^+$ and $A^{\dagger}$, obviously, obey the conditions of Lemmas~\ref{lem.cauchy0} and~\ref{lem.tech}.

Notice that equations~\eqref{eq.conse1d} can be rewritten as follows
\begin{gather}
\dot x_i=\sum\limits_{j\in\V^+}a_{ij}(t)(\hat x_j^i(t)-x_i(t))+f_i(t),\,i\in\V^+;\label{eq.subsys1+}\\
\dot x_i=\sum\limits_{j\in\V^{\dagger}}a_{ij}(t)(\hat x_j^i(t)-x_i(t))+f_i(t),\,i\in\V^{\dagger}\label{eq.subsys2+},\\
f_i(t)\dfb
\begin{cases}
\sum_{j\in\V^{\dagger}}a_{ij}(t)[\hat x^i_j(t)-x_i(t)],\quad i\in\V^+\\
\sum_{j\in\V^+}a_{ij}(t)[\hat x^i_j(t)-x_i(t)],\quad i\in\V^{\dagger}.
\end{cases}\label{eq.fi}
\end{gather}

Introducing the subvectors $x^+=(x_i)_{i\in\V^+}$ and  $x^{\dagger}=(x_i)_{i\not\in\V^{\dagger}}$ and corresponding vectors $f^+,f^{\dagger}$ and using
Lemma~\ref{lem.cauchy0} (where $t_*=\tau_m,U=U^+$), we have
\ben%\label{eq.aux4a}
\begin{gathered}
x^+(t)\leq U^+(t,\tau_m)x^+(\tau_m)+
\La(\tau_m)(\ones_{\V^+}-U^+(t,\tau_m)\ones_{\V^+})\\+\int_{\tau_m}^{t}U^+(t,\xi)f^+(\xi)d\xi\quad\forall t\geq\tau_m.
\end{gathered}
\een
Recalling that $x^+(\tau_m)\leq(\La(\tau_m)-\delta)\ones_{\V^+}$, one arrives at
\be\label{eq.aux4a}
\begin{gathered}
x^+(t)\leq (\La(\tau_m)-\delta)\ones_{\V^+}+\int_{\tau_m}^{t}U^+(t,\xi)f^+(\xi)d\xi.
\end{gathered}
\ee
Recalling that $x^{\dagger}(\tau_m)\leq\La(\tau_m)\ones_{\V^{\dagger}}$ and using Lemma~\ref{lem.cauchy0} (with $\tau_m,U^{\dagger}$ instead of $t_*,U$), one also proves that
\be\label{eq.aux4b}
x^{\dagger}(t)\leq \La(\tau_m)\ones_{\V^{\dagger}}+\int_{\tau_m}^{t}U^{\dagger}(t,\xi)f^{\dagger}(\xi)d\xi.
\ee

\textbf{Step 2.} Using the connectivity of $\g_{\infty}$, we will show that $\max x^+(t)\geq[\La(\tau_m)-2\delta/3]$ for some $t\geq\tau_m$. Indeed, assume on the contrary
$x^+(t)\leq[\La(\tau_m)-2\delta/3]\ones_{\V^+}$ at every instant $ t\geq\tau_m$. By assumption, $z_{r+1}(t)>\La(\tau_m)-\delta/3$ for all $t\geq\tau_m-\bar h$. Therefore,
$x_j(t)>\La(\tau_m)-\delta/3\,\forall j\in\V^{\dagger}$ and
\[
f_j(t)\overset{\eqref{eq.fi}}{\leq} -\frac{\delta}{3}\sum_{i\in\V^+}a_{ij}(t)\leq 0\quad\forall j\in\V^{\dagger}\;\;\forall t\geq\tau_m-\bar h.
\]
Recalling that the persistent graph is connected and $a_{ij}
\not\in L_1[0,\infty)$ for some pair $i\in\V^+,j\in\V^{\dagger}$, we have $\int_{\tau_m}^{\infty}f_j(t)dt=-\infty$.
Applying Lemma~\ref{lem.tech} to matrix $U^{\dagger}$, one shows, due to~\eqref{eq.aux4b}, that
$
x_j(t)\xrightarrow[t\to\infty]{}-\infty,
$
which contradicts to the solution's boundedness.

\textbf{Step 3.} Let $t'\geq\tau_s$ be the first instant such that $x_i(t')=\La(\tau_s)-2\delta/3$ for some $i\in\V^+$; find $q$ such that
$t_{q-1}<t'\leq t_q$. Recalling the definition of $(t_p)$-boundedness and using~\eqref{eq.simple-est}
\ben%\label{eq.aux2}
\begin{gathered}
\int_{t'}^{t_q}\alpha_i(s)ds\overset{\eqref{eq.ell}}{\leq}(n-1)\ell\overset{\eqref{eq.simple-est}}{\Longrightarrow}
x_i(t)\leq (1-\theta_0)\La(t')+\theta_0x(t')\\ \forall t\in[t',t_q],\quad \theta_0\dfb e^{-(n-1)\ell}\,\forall i\in\V^+,
%$x(t_p)\leq \psi x(t')+(1-\psi)\La(t')\ones$,
\end{gathered}
\een
Therefore,  for all $i\in\V^+$ and $t\in [\tau_m,t_q]$ one has
\be\label{eq.aux+}
x_i(t)\leq \La(\tau_m)-\delta_1,\quad\delta_1\dfb\frac{2\delta}{3}\theta_0.
\ee
On the other hand, by construction $x_i(t')-x_i(\tau_m)\geq\delta/3$. Using~\eqref{eq.aux4a} and recalling that $U(t,\xi)$ is substochastic, one has
\[
\sum_{l\in\V^+}\int_{\tau_m}^{t'}f_i(t)dt=\sum_{\substack{l\in\V^+,\\j\in\V^{\dagger}}}\int_{\tau_m}^{t'}a_{lj}(t)[\hat x^l_j(t)-x_l(t)]dt\geq\frac{\delta}{3}.
\]
Since $\hat x^l_j(t)-x_l(t)\le D_0\dfb\La(t_*)-\la(t_*)$ due to Lemma~\ref{prop.bound},
\[
\sum_{\substack{l\in\V^+,\\j\in\V^{\dagger}}}\int_{\tau_m}^{t'}a_{lj}(t)dt\geq c_1\dfb\frac{\delta}{3D_0}.
\]
Therefore, in view of the NITS condition\footnote{Recall that $\tau_m=t_{p(m)}$ is a member of the sequence $(t_p)$} and $t_q\geq t'$,
\be\label{eq.auxaux}
\sum_{\substack{l\in\V^+,\\j\in\V^{\dagger}}}\int_{\tau_m}^{t_q}a_{jl}(t)dt\geq c_2\dfb K^{-1}c_1.
\ee

Informally speaking, we have proved that the influence between $\V^+$ and $\V^{\dagger}$ during $[\tau_m,t_q]$ is strong enough. A final step is needed to prove~\eqref{eq.aux++}.

\textbf{Step 4.} Using~\eqref{eq.auxaux} and Lemma~\ref{lem.tech}, inequality~\eqref{eq.aux4b} entails that for
some $j\in\V^{\dagger}$ and $t''\in [\tau_m,t_q]$ we have $x_{j}(t'')\leq \La(\tau_m)-\delta_2$, where, by definition
\be\label{eq.const-aux+}
\begin{gathered}
\delta_2\dfb\frac{\varrho^{\dagger}c_2\delta_1}{\varrho^{\dagger}c_2+n-r},
\end{gathered}
\ee
$c_2$ is defined in~\eqref{eq.auxaux} and $\rho^{\dagger}$ is the constant from Lemma~\ref{lem.tech}, corresponding to $U^{\dagger}$.

Indeed, assume that $x_{j}(t)> \La(\tau_m)-\delta_2$ for all $j\in\V^{\dagger}$ and $t\in [\tau_m, t_q]$. In view of~\eqref{eq.aux-assum}, the latter inequality holds also for
$t\in [\tau_m-\bar h, \tau_m]$ (recall that $\delta_2<\delta/3$). Therefore, we have $\hat x^j_i(t)-x_j(t)>\delta_1-\delta_2$ for $i\in\V^+,j\in\V^{\dagger},t\in [\tau_m,t_q]$. Thus
\[
\sum_{j\in\V^{\dagger}}\int_{t_p}^{t_q}f_j(t)dt\overset{\eqref{eq.auxaux}}{<} -c_2(\delta_1-\delta_2)<0.
\]
Applying Lemma~\ref{lem.tech} to matrix $U^{\dagger}$, one has
\[
\begin{split}
(n-r)(\La(\tau_m)-\delta_2)\leq
\ones_{\V^{\dagger}}^{\top}x^{\dagger}(t_q)\overset{\eqref{eq.aux4b}}{\leq} (n-r)\La(\tau_m)+\\+
\int_{\tau_m}^{t_q}\ones_{\V^{\dagger}}^{\top}U^{\dagger}(t_q,t)f^{\dagger}(t)dt<\\
<(n-r)\La(\tau_m)-\rho^{\dagger}c_2(\delta_1-\delta_2).
\end{split}
\]
The latter inequality, obviously, contradicts~\eqref{eq.const-aux+}.

Now, let $\tilde p(m)\leq q$ be such an index that $t''\in(t_{\tilde p(m)-1},t_{\tilde p(m)}]$. Using~\eqref{eq.ell},
the inequality $x_{j}(t'')\leq \La(\tau_m)-\delta_2$ entails, similarly to~\eqref{eq.aux+}, that
\[
x_j(t_{\tilde p(m)})\leq\La(\tau_m)-\tilde\delta,\quad\tilde\delta\dfb\theta_0\delta_2.
\]
The latter inequality, in combination with~\eqref{eq.aux+}, obviously entails~\eqref{eq.aux++}, where $\tilde\tau_m=t_{\tilde p(m)}$, which contradicts to the induction hypothesis. 

The contradiction shows that~\eqref{eq.z-consensus1} is valid, which finishes the proof in the case of connected graph
$\g_{\infty}$.

\begin{rem}
In fact, the proof of Theorem~\ref{thm.symm-c} remains valid (modulo minor modifications) for any bounded solution of~\eqref{eq.conse1d-ineq}; the sketch of the proof for inequalities can be found in~\cite{ProCala_CDC:2020}.
\end{rem}

\subsubsection{The case of \textbf{disconnected} persistent graph $\g_{\infty}$}

In the general situation, convergence of solutions and the consensus in each connected component is immediate from Lemma~\ref{lem-l1-rob}.
To reuse some equations from the previous section, denote the set of nodes belonging to some connected component by
 $\V^{\dagger}\subseteq\V$ and let $\V^+=\V\setminus\V^{\dagger}$. In the previous subsection, we have already proved that the ``reduced'' algorithm~\eqref{eq.subsys2}
 establishes consensus among the agents from $\V^{\dagger}$. Notice also that, by definition of the persistent graph, $a_{ij}\in L_1$ for all $i\in\V^+,j\in\V^{\dagger}$,
 therefore, all $f_i$ in~\eqref{eq.fi} are also $L_1$-summable (recall that the solution is bounded). Lemma~\ref{lem-l1-rob} now entails that equations~\eqref{eq.subsys2+} also imply consensus between
 each two agents from $\V^{\dagger}$, which finishes the proof of Theorem~\ref{thm.symm-c}.$\blacksquare$

 \subsection{Proofs of Theorems~\ref{thm.repeated-disc} and~\ref{thm.symm-d}}\label{subsec.proof-discrete}

 Theorems~\ref{thm.repeated-disc} and~\ref{thm.symm-d} reduce to Theorems~\ref{thm.repeated-cont}
 and~\ref{thm.symm-c} respectively by applying the following trick that if often used in analysis of sampled-time systems~\cite{Fridman:88,SelivanovFridman_Survey:2019}
 and reduces a discrete-time system to a \emph{continuous-time} delay system with time-varying ``sawtooth'' delays.
 Consider the discrete-time algorithm\footnote{It is convenient to change notation here, denoting the discrete
time by $k$, the weight matrix by $B(k)$ and the solution by $z(k)$.}
 \be\label{eq.eq-del-disc}
 \begin{gathered}
z_i(k+1)=z_i(k)+\sum_{j\ne i}b_{ij}(k)[\hat z_j^i(k)-z_i(k)],\quad \forall i\in\V,\\
\hat z^i_j(k)\dfb z_j(k-h_{ij}^0(k)).
\end{gathered}
\ee

\begin{lem}\label{lem.b-to-a}
Consider the sequence of stochastic matrices $B(k)$ with positive diagonal entries $b_{ii}(k)>0$.
Consider a piecewise-constant function matrix $A(t)$, defined by
\be\label{eq.b-to-a}
\begin{aligned}
a_{ij}(t)&=\frac{-b_{ij}(k)\ln b_{ii}(k)}{1-b_{ii}(k)}\quad&\forall t\in[k,k+1)\,\forall i\ne j,\\
a_{ii}(t)&=0&\forall t\geq 0\,\forall i.
\end{aligned}
\ee
Consider also ``sawtooth'' delay functions
\[
h_{ij}(t)=t-k+h_{ij}(k)\quad\forall i\ne j\,\forall t\in[k,k+1).
\]
Then, for any solution to~\eqref{eq.eq-del-disc},
the function $x(k)$ defined on each interval $[k,k+1)$ as a solution to the Cauchy problem
\be\label{eq.z-to-x1}
%\begin{gathered}
\dot{x}_i(t)=\sum_{j\ne i}a_{ij}(t)[\hat z_j^i(k)-x_i(t)],\quad x_i(k)=z_i(k),
%\end{gathered}
\ee
is absolutely continuous and obeys~\eqref{eq.conse1d}. If $B$ is strongly aperiodic, then the
matrix $A(\cdot)$ is bounded (and, in particular, Assumption~\ref{asm.int-bound} holds).
If $B(\cdot)$ fulfils the AQSC or NITS condition, the respective condition is also satisfied by $A(\cdot)$.
\end{lem}
\begin{IEEEproof}
A straightforward computation shows that $\alpha_{i}(t)\dfb\sum_{j\ne i}{a_{ij}(t)}\equiv -\ln b_{ii}(k)\,\forall t\in[k,k+1)$ and
\[
\begin{gathered}
\sum_{j\ne i}a_{ij}(k)\hat z_j^i(k)=\frac{-\ln b_{ii}(k)}{1-b_{ii}(k)}[z_i(k+1)-b_{ii}(k)z_i(k)]=\\
=\frac{\alpha_i(k)}{1-e^{-\alpha_i(k)}}[z_i(k+1)-b_{ii}(k)z_i(k)].
\end{gathered}
\]
Solving equation~\eqref{eq.z-to-x1}, one obtains
\[
\begin{split}
\lim_{\substack{t\to (k+1)\\t<k+1}}x_i(t)=e^{-\alpha_i(k)}z_i(k)+
\sum_{j\ne i}\frac{1-e^{-\alpha_i(k)}}{\alpha_i(k)}a_{ij}(k)\hat z_j^i(k)=\\=\sum_{j\ne i}b_{ij}z_j(k-h_{ij}^0(k))=z_i(k+1)=x_i(k+1).
\end{split}
\]
Hence, $x(t)$ is a continuous function at $t=k=1,\ldots$, being thus also absolutely continuous.
By noticing that $z_j(k-h_{ij}^0(k))=x_j(t-h_{ij}(t))$,~\eqref{eq.z-to-x1} entails~\eqref{eq.conse1d}.
The remaining two statements are immediate from~\eqref{eq.b-to-a} and Definitions~\ref{def.tp-bound}, \ref{def.aqsc} and
\ref{def.symm}: obviously, function $A(\cdot)$ inherits the properties of $(t_p)$-boundedness, AQSC and NITS
from $B(\cdot)$, provided that diagonal entries $b_{ii}(k)$ are uniformly positive.
\end{IEEEproof}

\textbf{Proof of Theorem~\ref{thm.repeated-disc}.} Theorem~\ref{thm.repeated-disc} now trivially follows from Theorem~\ref{thm.repeated-cont},
because algorithm~\eqref{eq.conse0d}, under the assumptions of strong periodicity and AQSC, reduces to algorithm~\eqref{eq.conse1d}
that obeys Assumption~\ref{asm.int-bound} and the AQSC condition. $\blacksquare$

\textbf{Proof of Theorem~\ref{thm.symm-d}.} Similarly, Theorem~\ref{thm.symm-d} follows from Theorem~\ref{thm.symm-c},
because algorithm~\eqref{eq.conse0d}, under the assumptions of strong periodicity and NITS,
reduces to algorithm~\eqref{eq.conse1d} that obeys Assumption~\ref{asm.int-bound} and the NITS condition. $\blacksquare$

\subsection{Proofs of Theorems~\ref{thm.nlin-c},~\ref{thm.nlin-d}}

Consider an arbitrary solution $x(t)$ to~\eqref{eq.conse1nd} or~\eqref{eq.conse0nd}. For this special solution, we introduce
the matrix $\hat A$ as follows
\[
\hat a_{ij}(t)=a_{ij}(t)\psi_{ij}(\hat x^i_j(t),x_i(t)),\,\forall i\ne j
\]
and $\hat a_{ii}=1-\sum_{j\ne i}\hat a_{ij}\,\forall i\in\V$ in the discrete-time case (in the continuous-time case,
$\hat a_{ii}=0$ without loss of generality).
Due to~\eqref{eq.phi-struct}, the solution $x(t)$ obeys respectively~\eqref{eq.conse1nd} or~\eqref{eq.conse0nd},
where $a_{ij}$ has to be replaced by $\hat a_{ij}$. Therefore, Lemma~\ref{prop.bound} ensures that $x(t)$ is bounded; in view of
~\eqref{eq.psi-ineq}, $\psi_{ij}(\hat x^i_j(t),x_i(t))$ is uniformly positive and bounded:
\[
0<m\leq \psi_{ij}(\hat x^i_j(t),x_i(t))\leq M<\infty\quad\forall t.
\]
Also, in the discrete-time case we have $\hat a_{ii}\geq a_{ii}\,\forall i$. For this reason, it can be
shown that matrix $A$ satisfies the condition of
Theorems~\ref{thm.repeated-cont}-\ref{thm.symm-d} if and only if $\hat A$ abides by these conditions.$\blacksquare$

\subsection{Proofs of Lemmas~\ref{lem.convex},~\ref{Lem.contain-stab} and~\ref{Lem.aggreg-stab}}

In this subsection, we prove Lemmas~\ref{lem.convex},~\ref{Lem.contain-stab} and~\ref{Lem.aggreg-stab} that are closely related.

\subsubsection{Proof of Lemma~\ref{lem.convex}}

Lemma~\ref{lem.convex} is, in fact, a multidimensional extension of Lemma~\ref{prop.bound} and can be derived from it.
Assume on the contrary that $x_i(t)\not\in\mathcal{D}$ for some $i\in\V$ and $t>t_*$. Due to the Hahn-Banach theorem,
the vector $f\in\r^m$ exists such that $\xi_*\dfb f^{\top}x_i(t)>f^{\top}v\,\forall v\in\mathcal{D}$. At the same time, the values
$\xi_i(t)=f^{\top}x_i(t)$ obey the \emph{scalar} equations~\eqref{eq.conse0d} or~\eqref{eq.conse1d} and
$\xi_i(t)<\xi_*$ when $t-\bar h\leq t\leq t$ (and thus $\La(t_*)<\xi_*$).
Applying Lemma~\ref{prop.bound}, one shows that $\xi_i(t)<\xi_*\,\forall t\geq t_*$, which contradicts to our assumption.
Lemma~\ref{lem.convex} is proved. $\blacksquare$

\subsubsection{Proof of Lemma~\ref{Lem.aggreg-stab}}

The proof of ``only if'' part is straightforward, choosing $\Omega=\{0\}$ and $\omega_i\equiv 0$.

To prove the ``if'' part, we show first that the target set $\Omega^{\V}$ is forward invariant, that is, if
$x_i(t)\in\Omega\,\forall i\in\V\,\forall t\leq t_*$, then $x_i(t)\in\Omega\,\forall i\in\V\,\forall t\geq t_*$.
Repeating the arguments from~\cite[Section~4.1]{ProCao:2017}, it can be shown that the functions
$\xi_i(t)=\dist(x_i,\Omega)$ obey the differential averaging inequality~\eqref{eq.conse1d-ineq}.
As shown in~\cite[Proposition~1]{ProCala_CDC:2020}, the function $\La(t)\dfb\sup_{s\in[t-\bar h,t]}\max_i\xi_i(s)$ is therefore
non-increasing, in particular, if $\xi_i=0\,\forall t\leq t_*$, then $\xi_i\equiv 0$.

Notice now that stability of~\eqref{eq.avg-cont0-damp} means that the difference between any two solutions $(x_i)$ and $(\tilde x_i)$ to~\eqref{eq.avg-cont0-contain} vanishes
as $t\to\infty$. Considering a solution $\tilde x$ such that $\tilde x_i(t)\in\Omega\,\forall t\leq t_*$, we have
$\tilde x_i(t)\in\Omega\,\forall t\geq t_*$, and hence any other
solution $(x_i(t))$ converges~\eqref{eq.aggreg} to set $\Omega$. $\blacksquare$

\subsubsection{Proof of Lemma~\ref{Lem.contain-stab}} To prove the ``only if'' part is suffices to consider the special
case where all leaders are located at the same position $x_k=0\,\forall k\in\V$; in this situation,~\eqref{eq.avg-cont0-damp}
coincides with~\eqref{eq.avg-cont0-contain}, and containment~\eqref{eq.containment} means that all solutions converge to $0$.

The ``if'' part is immediate from Lemma~\ref{Lem.aggreg-stab} by noticing that~\eqref{eq.avg-cont0-contain} is a special case
of~\eqref{eq.avg-cont0-agreg}, where $\Omega=\mathcal{C}_L$ and
\[
\omega_i(t)=\frac{1}{d_i(t)}\sum_{k\in\V_L}b_{ik}x_k\
\]
(if $d_i(t)=0$, then $\omega_i(t)$ is an arbitrary element of $\Omega$). $\blacksquare$

\section{Conclusion}\label{sec.concl}

First-order consensus protocols, based on the idea of iterative averaging (discrete time) or the Laplacian flow (continuous-time),
are prototypic distributed algorithms that arise in many problems of multi-agent coordination. In spite of substantial progress
in their analysis, some fundamental problems still remain open, in particular, there are no necessary and sufficient criteria
of the algorithm's convergence in the situation where the interaction graph is time-varying.
The latter problem becomes especially complicated in the case where communication among the agents is delayed.
The most typical assumption under which the delay robustness of consensus is proved the uniform quasi-strong connectivity (UQSC).

In this paper, we demonstrate that the UQSC assumption is not necessary and can be substantially relaxed. In the case of
a general directed graph, it can be relaxed to the condition we refer to as the aperiodic quasi-strong connectivity (AQSC),
which, as we show, is in fact very close to the well-known necessary consensus condition (integral connectivity).
In the case of reciprocal (e.g. undirected or type-symmetric), the necessary condition in fact also appears to be sufficient
(provided that the delays are bounded, and some other technical assumptions hold). We show that the discrete-time and
continuous-time algorithms in the delayed case can be studied simultaneously, because the continuous-time algorithm in fact
includes discrete-time averaging dynamics as a special case.

We also show that some algorithms for coordination of mobile agents such as containment control and target aggregation,
in fact, also reduce to consensus algorithms being thus covered by the general theory we develop. The latter fact has not been
realized in the existing literature, and containment control protocols are usually studied separately from consensus algorithms.
Our approach not only simplifies their analysis, but allows to generalize a number of results available in the literature.

\appendices
\section{Discussion on Assumption~\ref{asm.int-bound}}\label{sec.app-a}

In this appendix, we demonstrate that Assumption~\ref{asm.int-bound}
in fact cannot be omitted even for the case of two agents.
Consider two agents whose values evolve in accordance with
\be\label{eq.aux31}
\begin{gathered}
\dot x_0(t)=a(t)(x_1(t-\tau(t))-x_0(t)),\\
\dot x_1(t)=a(t)(x_0(t-\tau(t))-x_1(t)).
\end{gathered}
\ee
Here $a(t)\geq 0$ is a locally $L_1$-summable function, however $a\not\in L_1[0,\infty)$.
As shown in Lemma~\ref{prop.intermittent}, the AQSC condition is satisfied for some sequence $(t_p)$;
the NITS condition is also fulfilled with the same sequence.
All conditions of Theorems~\ref{thm.repeated-cont} and~\ref{thm.symm-c} thus holds except for Assumption~\ref{asm.int-bound}.

We now show that the function $a(\cdot)$ and delay $\tau(\cdot)$ can be chosen in such a way that the system~\eqref{eq.aux31} fails to establish consensus.
We define the delay in the same way as in Subsect.~\ref{subsec.proof-discrete}, that is, $t-\tau(t)=k$ and
\be\label{eq.aux31+}
\dot x_i(t)=a(t)(x_{1-i}(k)-x_i(t)),\; i=0,1,\; \forall t\in[k,k+1).
\ee
Denoting $a_k\dfb\exp(-\int_{k}^{k+1}a(t)\,dt)$, one has
\[
x_{1-i}(k)-x_i(k+1)=a_k(x_{1-i}(k)-x_i(k)),
\]
which leads to a discrete-time consensus algorithm as follows
\be\label{eq.aux31++}
%\begin{gathered}
x_i(k+1)=(1-a_k)x_{1-i}(k)+a_kx_i(k),\quad i=0,1,
%\end{gathered}
\ee
which, as easily shown, ensures that $|x_1(k+1)-x_0(k+1)|=|1-2a_k|\,|x_1(k)-x_0(k)|$.
Obviously, if $a_k<1/2$ and $\sum_k a_k<\infty$, then $x_1(k)-x_0(k)\nrightarrow 0$ as $k\to\infty$, so the algorithm fails to reach consensus (notice that the discrete-time dynamics~\eqref{eq.aux31++} fails
to satisfy the strong aperiodic condition).

\bibliographystyle{IEEETran}
\bibliography{consensus,social}

\end{document}